\documentclass[11pt,a4paper]{amsart}

\usepackage{amssymb}
\usepackage{microtype}

\usepackage{xcolor}
\definecolor{linkblue}{RGB}{38,68,105}
\definecolor{citegreen}{RGB}{45,88,68}
\usepackage[colorlinks=true,linkcolor=linkblue,citecolor=citegreen,urlcolor=linkblue,pdfborder={0 0 0}]{hyperref}

\newtheorem{theorem}{Theorem}[section]
\newtheorem{proposition}[theorem]{Proposition}
\newtheorem{corollary}[theorem]{Corollary}
\newtheorem{lemma}[theorem]{Lemma}
\theoremstyle{definition}
\newtheorem{definition}[theorem]{Definition}
\theoremstyle{remark}

\numberwithin{equation}{section}

\newcommand{\C}{\mathbb C}
\newcommand{\R}{\mathbb R}
\newcommand{\N}{\mathbb N}
\newcommand{\Z}{\mathbb Z}
\newcommand{\PP}{\mathbb P}
\newcommand{\HH}{\mathbb H}
\newcommand{\X}{\mathcal X}
\newcommand{\A}{\mathcal A}
\newcommand{\K}{\mathbb K}
\newcommand{\U}{\mathcal U}
\newcommand{\dbar}{\bar\partial}
\DeclareMathOperator{\SL}{SL}
\DeclareMathOperator{\PSL}{PSL}
\DeclareMathOperator{\Hol}{Hol}
\DeclareMathOperator{\Hom}{Hom}
\DeclareMathOperator{\End}{End}
\DeclareMathOperator{\tr}{tr}
\DeclareMathOperator{\Hopf}{Hopf}
\DeclareMathOperator{\diag}{diag}

\title[Boundary fibers of the holonomy map]{Boundary fibers of the holonomy map for complex projective structures}

\author{Shiyu Cao}
\address{Chern Institute of Mathematics and LPMC, Nankai University, Tianjin 300071, China}
\email{shiyucao@126.com}

\author{Sheng Tan}
\address{School of Mathematical Sciences, Capital Normal University, Beijing 100048, China}
\email{tansheng2018@outlook.com}

\date{}

\begin{document}

\begin{abstract}
We study the Morgan--Shalen boundary of holonomy representations of complex projective structures on a compact Riemann surface of genus at least two. We prove that every direction defined by a holomorphic quadratic differential with a multiple zero admits distinct projective marked length limits. Together with Dumas's uniqueness theorem for simple zeros, this gives a complete characterization of directions with unique limits. We also establish that no nonzero abelian length function can occur in the holonomy boundary of any such fixed surface, regardless of the limiting quadratic differential. Finally, we characterize all directional boundary fibers by a necessary-and-sufficient analytic condition involving valued representations and the holomorphic jet bundle underlying an oper.
\end{abstract}

\maketitle
\tableofcontents


\section{Introduction}\label{sec:introduction}

A complex projective structure on a compact Riemann surface determines a holomorphic quadratic differential through its Schwarzian derivative and a representation of the fundamental group through its holonomy \cite{GKM00,Dum17}. The behavior of these objects at infinity is related to measured foliations and actions on real trees \cite{MS84,DDW00,Dum17}. Dumas proved that when the limiting quadratic differential has only simple zeros, its horizontal measured foliation determines a unique projective marked length limit \cite[Theorem~B]{Dum17}. When the differential has a multiple zero, however, a straight map from its dual tree to a limiting tree may identify branches \cite[Theorem~A]{Dum17}. 
These results raise two questions about the boundary fibers of the holonomy map: whether distinct limits can arise in directions with multiple zeros and whether nonzero abelian length functions can occur.


\subsection{Main results}

Let $X$ be a compact connected Riemann surface of genus $g\ge2$, and let $K_X$ denote its canonical line bundle. Write $\Gamma=\pi_1(X)$ and $Q(X)=H^0(X,K_X^{\otimes2})$, equipped with the $L^1$ norm $\|\varphi\|_1=\int_X|\varphi|$. 
A choice of reference projective structure identifies the space of projective structures on $X$ with $Q(X)$. A theta characteristic, that is, a square root of $K_X$, determines a coherent lift of projective holonomy to $\SL_2(\C)$. We denote the resulting map to the $\SL_2(\C)$ character variety by
\[
\Hol_X:Q(X)\longrightarrow\X(\Gamma).
\]

For a representation $\rho:\Gamma\to\SL_2(\C)$, let $\ell_\rho(\gamma)$ be the translation length of $\rho(\gamma)$ on hyperbolic three-space of curvature $-1$. Boundary points in the Morgan--Shalen compactification are represented by nonzero marked length functions $\ell:\Gamma\to\R_{\ge0}$, considered up to multiplication by a positive constant \cite{MS84,DDW00}. Each such length function arises from an isometric action of $\Gamma$ on a real tree.

\begin{definition}\label{def:intro-fiber}
For $q\in Q(X)$ with $\|q\|_1=1$, the \emph{directional holonomy fiber} $\mathcal L_X(q)$ consists of the projective length functions $[\ell]$ for which there is a sequence $\varphi_j\in Q(X)$ satisfying
\[
\|\varphi_j\|_1\longrightarrow\infty,\quad \frac{\varphi_j}{\|\varphi_j\|_1}\longrightarrow q,\quad \Hol_X(\varphi_j)\longrightarrow[\ell]
\]
in the Morgan--Shalen compactification.
\end{definition}

The direction in this definition is a positive real ray. In particular, the horizontal foliation generally changes when $q$ is multiplied by a complex number of nonzero argument. For a nonzero differential $q$, let $T_q$ be the tree dual to its horizontal measured foliation on the universal cover, with transverse measure $|\operatorname{Im}\sqrt q|$.

Dumas proved that every holonomy limit in the direction of $q$ is represented by a tree admitting a surjective $\Gamma$-equivariant straight map from $T_q$ \cite[Theorem~A]{Dum17}. Such a map preserves the transverse metric along the images of nonsingular flat geodesic segments. When all zeros of $q$ are simple, he further proved that the projective marked length limit is unique and equals $[\ell_{T_q}]$ \cite[Theorem~B]{Dum17}. When $q$ has a multiple zero, however, a straight map may identify branches of $T_q$, raising the question of whether distinct holonomy limits can arise in the same direction. Dumas also asked whether a nonzero abelian length function can occur as a holonomy limit on a given Riemann surface \cite{Dum17}. We prove that every direction corresponding to a differential with a multiple zero admits distinct limits and that no nonzero abelian length function can occur.

\begin{theorem}\label{thm:intro-nonunique}
Let $X$ be a compact Riemann surface of genus at least two, and let $q\in Q(X)\setminus\{0\}$ be a nonzero holomorphic quadratic differential on $X$ with a multiple zero. Then there exist $\gamma,\delta\in\Gamma$ and constants $H,h>0$ such that, for every $\lambda\in[0,\sqrt2H]$, there are $t_j\to\infty$ and $\varphi_j\in Q(X)$ with
\[
\frac{\varphi_j}{t_j^2}\longrightarrow q,\quad \frac{\ell_{\Hol_X(\varphi_j)}(\beta)}{t_j}\longrightarrow\ell(\beta)\quad(\beta\in\Gamma),
\]
where $\ell(\gamma)=\lambda$ and $\ell(\delta)=\sqrt2h$. 
Consequently, for every $r\in[0,H/h]$, there exists $[\ell]\in\mathcal L_X(q/\|q\|_1)$ such that $\ell(\delta)>0$ and $\ell(\gamma)/\ell(\delta)=r$.
\end{theorem}

The constructed length functions take the same positive value on $\delta$ but different values on $\gamma$, so their projective classes are distinct. Combining this with Dumas's uniqueness theorem for quadratic differentials whose zeros are all simple \cite[Theorem~B]{Dum17}, we obtain the following characterization.

\begin{corollary}\label{cor:intro-dichotomy}
For every nonzero $q\in Q(X)$,
\[
\bigl|\mathcal L_X(q/\|q\|_1)\bigr|=
\begin{cases}
1,&\text{all zeros of $q$ are simple},\\
2^{\aleph_0},&\text{$q$ has a multiple zero}.
\end{cases}
\]
In the first case the fiber is $\{[\ell_{T_q}]\}$.
\end{corollary}

Our next main result concerns whether nonzero abelian length functions can arise as holonomy limits on $X$.

\begin{definition}\label{def:intro-abelian}
A nonzero length function $\ell$ on $\Gamma$ is \emph{abelian} if it has the form $\ell(\gamma)=|c(\gamma)|$ for a nonzero homomorphism $c:\Gamma\to\R$.
\end{definition}

An abelian length function arises from an action of $\Gamma$ on $\R$ by translations. The corresponding equivariant harmonic function on the universal cover of $X$ has a differential whose $(1,0)$ part is holomorphic \cite{DDW00}. Dumas proved that if an abelian length function occurs as a holonomy limit, then the limiting quadratic differential lies on the positive ray of $\omega^2$, where $\omega$ is a holomorphic one-form whose imaginary periods are proportional to $c$ \cite[Theorem~C]{Dum17}. The following theorem excludes abelian limits from the image of the holonomy map.

\begin{theorem}\label{thm:intro-exclusion}
For every compact Riemann surface $X$ of genus at least two, no nonzero abelian length function represents an accumulation point of $\Hol_X(Q(X))$ in the Morgan--Shalen boundary of the character variety.
\end{theorem}

Consequently, no nonzero class in $\PP H^1(X,\R)$ can arise from a holonomy limit. Applying the theorem to finite unramified covers of $X$, we further show that every nontrivial minimal tree action arising from a holonomy limit preserves neither an end nor a line. In particular, if $q=\omega^2$ for a nonzero holomorphic one-form $\omega$, Theorem~\ref{thm:intro-nonunique} gives continuum many projective length functions in $\mathcal L_X(q/\|q\|_1)$, but none is represented by the abelian length function $\gamma\mapsto|\operatorname{Im}\int_\gamma\omega|$.


\subsection{Related work}

The relation between representation limits and tree actions originates in the work of Culler--Shalen and Morgan--Shalen \cite{CS83,MS84}. 
Marked-length rigidity for irreducible minimal actions was established by Culler--Morgan \cite{CM87}. Harmonic maps connect this boundary to quadratic differentials: Wolf studied maps from surfaces to trees and the realization of measured foliations \cite{Wol95,Wol96}, and Farb--Wolf analyzed the corresponding folding and splitting phenomena \cite{FW01}. 
Daskalopoulos--Dostoglou--Wentworth constructed harmonic limits of degenerating representations and proved that the normalized Hopf differential is determined by the projective length function \cite{DDW98,DDW00}. 
He--Wentworth--Zhang further studied the relation between $\Z/2$ harmonic one-forms and actions on real trees \cite{HWZ26}. These results concern the geometry of possible limiting actions.
The holonomy representations of projective structures on $X$ are subject to an additional global holomorphic constraint.

Gallo--Kapovich--Marden proved that every nonelementary representation of a closed surface group into $\PSL_2(\C)$ that lifts to $\SL_2(\C)$ is the holonomy of a projective structure on some Riemann surface \cite{GKM00}. 
When the complex structure is prescribed, the holomorphic bundle associated with the flat connection must be isomorphic to the oper jet bundle. 
Dumas--Kent proved that Bers slices are Zariski dense \cite{DK09}, while Baba studied intersections of holonomy varieties associated with different conformal structures \cite{Bab23,Bab25}. 
The finiteness theorem for Fuchsian equations with regular singularities in \cite{Tan26} also uses harmonic periods to compare limiting quadratic differentials in the abelian case. 
For compact projective structures, our argument excludes abelian limits by constructing approximately holomorphic sections of a line bundle of negative degree and applying a uniform Dolbeault estimate.

Ott--Swoboda--Wentworth--Wolf relate limiting configurations to harmonic maps and pleated surfaces and establish asymptotics for oper rays whose leading differentials have only simple zeros \cite{OSW$^+$24}. Their proof of Proposition~6.9 uses character correction, estimates for small analytic errors, and elliptic compactness on the oper bundle. To exclude abelian limits, we combine these methods with an algebraic argument. The vanishing of normalized lengths for commutators involving every positive power of a chosen group element forces the generators to act approximately by scalars on a common line, with errors decaying faster than any exponential rate. Estimates for transport of the scalar equation then apply across all zeros of the limiting differential.

For the proof of nonuniqueness, we use Sibuya's theory of canonical solutions and Stokes multipliers \cite{Sib75}, together with Shin's results on the linear dependence of solutions decaying in nonadjacent sectors \cite{Shi04}. The quadratic case is governed by the Weber equation \cite[Chapter~12]{OOL$^+$26}. WKB connection formulas and Stokes graphs are developed in \cite{IN14}, while Acosta studied how local asymptotics can be matched globally for cyclic opers \cite{Aco16}. Products of Stokes matrices also appear in the asymptotics of cubic and quartic Hitchin rays \cite{LTW26,Ma26}, and more general estimates for harmonic maps at high energy are given in \cite{SS26}. In our construction, we control the entire connection matrix uniformly even at a spectral resonance, where its leading scalar coefficient vanishes. We then use a global holomorphic parameter to prescribe the exact holonomy trace along a chosen closed curve.

Our final result characterizes which projective marked length functions arise in a given directional fiber. We establish a necessary and sufficient analytic criterion involving valued representations and the holomorphic bundle underlying an oper (Theorem~\ref{thm:criterion}). We use a coefficient field based on Robinson's theory of asymptotic numbers \cite{Rob66,TW04}, which is related to differential algebras developed for generalized functions \cite{Tod99,KSV05}. Oberguggenberger--Pilipovi\'c--Valmorin and Khalfallah--Kosarew have also studied holomorphic representatives in related settings \cite{OPV07,KK18}. The condition on the bundle draws on the classical theory of holomorphic connections and opers \cite{Ati57,Wen16,BHS23}. To recover the limiting quadratic differential, we use Dumas's results on straight maps \cite{Dum17} and the uniqueness of the Hopf differential for nonabelian tree actions \cite{DDW00}.

The outline of the paper is as follows. In Section~\ref{sec:preliminaries}, we recall the relevant facts about opers and length functions, together with the estimates needed later. Section~\ref{sec:exclusion} proves the exclusion of abelian limits and derives its consequences for finite unramified covers. We establish nonuniqueness for quadratic differentials with multiple zeros in Section~\ref{sec:nonuniqueness}. Finally, Section~\ref{sec:realization} gives an analytic criterion for determining which length functions arise as holonomy limits in a given direction.

\section{Opers and the holonomy boundary}\label{sec:preliminaries}

For the rest of the paper, let $X$ be a compact connected Riemann surface of genus $g\ge2$, with canonical line bundle $K_X$, and choose a base point $x_0\in X$. Write $\Gamma=\pi_1(X,x_0)$, and let $\widetilde X$ be the universal cover of $X$. We equip $Q(X)=H^0(X,K_X^{\otimes2})$ with the $L^1$ norm $\|\varphi\|_1=\int_X|\varphi|$ and choose a theta characteristic $L$ with $L^{\otimes2}=K_X$.


\subsection{The scalar equation and the jet bundle}\label{subsec:opers}

We begin by recalling the definition of a complex projective structure.

\begin{definition}\label{def:projective}
An \emph{complex projective structure} on $X$ is a maximal atlas of locally injective holomorphic charts to $\PP^1$ whose transition maps are restrictions of M\"obius transformations. Its developing map is a locally injective holomorphic map $f:\widetilde X\to\PP^1$ satisfying $f(\gamma x)=\rho(\gamma)f(x)$ for a homomorphism $\rho:\Gamma\to\PSL_2(\C)$, called its holonomy.
\end{definition}

\begin{definition}\label{def:projective-connection}
The \emph{Schwarzian derivative} of a locally injective holomorphic function $f$ in a coordinate $z$ is
\[
S(f;z)=\frac{f'''}{f'}-\frac{3}{2}\left(\frac{f''}{f'}\right)^2.
\]
A \emph{holomorphic projective connection} on $X$ is a collection of holomorphic functions $P_z$ in local coordinates $z$ satisfying
\[
P_z(z)=(w'(z))^2P_w(w(z))+S(w;z)
\]
under every holomorphic change of coordinate $w=w(z)$.
\end{definition}

Fix a projective structure with projective connection $P_0$. In a local coordinate $z$, the structure with Schwarzian parameter $\varphi\in Q(X)$ is represented by
\begin{equation}\label{eq:scalar}
y''+\frac{P_0+\varphi}{2}y=0.
\end{equation}
Here $P_0$ and $\varphi$ denote their local coefficients, and primes are derivatives with respect to $z$. The Schwarzian derivative of a ratio of independent solutions is $P_0+\varphi$. 
With $M=L^{-1}$, scalar solutions of \eqref{eq:scalar} are local holomorphic sections of $M$. The first jet of a local section records its value and first derivative. In a local coordinate $z$, the first jet of $y(z)(dz)^{-1/2}$ is represented by the pair $(y(z),y'(z))$. These jets form the rank-two holomorphic vector bundle
\begin{equation}\label{eq:jet}
E=J^1(M),\quad 0\longrightarrow L\longrightarrow E\longrightarrow M\longrightarrow0.
\end{equation}
The projection $E\to M$ takes a jet to its value, while its kernel is naturally identified with $K_X\otimes M\simeq L$. Since \eqref{eq:scalar} is a second-order equation, a local solution is uniquely determined by its first jet at any point. The determinant of $E$ is canonically trivial. The first-jet extension is nonsplit, since its extension class is the Atiyah class of $M$, which is nonzero because $\deg M=1-g\ne0$ \cite{Ati57}.
The determinant trivialization allows us to choose local bases of solutions with Wronskian one. Their analytic continuation gives an $\SL_2(\C)$ holonomy representation, well defined up to conjugacy, and hence the map 
\[
\Hol_X:Q(X)\to\X(\Gamma).
\]
The relation between projective holonomy and its lifts is discussed in \cite[Section~1.3]{GKM00}.

\begin{definition}\label{def:oper}
An $\SL_2(\C)$-\emph{oper} with theta characteristic $L$ is a tracefree holomorphic connection $\nabla$ on $E$, together with the line bundle $L\subset E$, such that its second fundamental form
\[
\beta_\nabla:L\longrightarrow(E/L)\otimes K_X
\]
is an isomorphism. The map $\beta_\nabla$ is obtained by applying $\nabla$ to a local section of $L$ and projecting modulo $L\otimes K_X$.
\end{definition}

The scalar equation \eqref{eq:scalar}, its associated holomorphic connection on $E$, and the corresponding oper are equivalent descriptions of the same projective structure \cite[Section~4.2]{BHS23}. Every holomorphic connection on a complex curve is flat. In particular, a tracefree holomorphic connection on $E$ defines an oper precisely when its second fundamental form is nonzero. For this jet bundle, the condition holds for every holomorphic connection: if the second fundamental form vanished, $L$ would inherit a holomorphic connection, contradicting $\deg L=g-1>0$. This is the standard degree argument in oper theory \cite[Section~4]{Wen16,BHS23}.

We next derive a consequence of the fact that the extension \eqref{eq:jet} does not split. Fix an oper $(E,L,\nabla_0)$ and let $B\in H^0(X,\End_0(E)\otimes K_X)$, where $\End_0(E)$ denotes the bundle of tracefree endomorphisms of $E$. The composite
\[
L\xrightarrow{B|_L}E\otimes K_X\longrightarrow M\otimes K_X\simeq L
\]
is multiplication by a constant. If this constant were nonzero, dividing $B|_L$ by it would split the extension obtained by tensoring \eqref{eq:jet} with $K_X$. Hence the composite vanishes, and
\begin{equation}\label{eq:unchanged-oper}
B(L)\subset L\otimes K_X,\quad \beta_{\nabla_0+B}=\beta_{\nabla_0}.
\end{equation}
Thus a global holomorphic correction of this form preserves the oper condition.


\subsection{Lengths, straight maps, and normalization}\label{subsec:lengths}

Equip hyperbolic three-space $\HH^3$ with its metric of curvature $-1$. The group $\PSL_2(\C)$ acts on $\HH^3$ by orientation-preserving isometries, and $\SL_2(\C)$ acts through its quotient $\PSL_2(\C)$. For $A\in\SL_2(\C)$, define its translation length by
\[
\ell(A)=\inf_{x\in\HH^3}d_{\HH^3}(x,Ax),
\]
where $Ax$ denotes the action of the isometry induced by $A$. For a representation $\rho:\Gamma\to\SL_2(\C)$, the associated length function on $\Gamma$ is defined by
\[
\ell_\rho(\gamma)=\ell(\rho(\gamma)),\quad \gamma\in\Gamma.
\]
If $r\ge1$ is the larger of the moduli of the eigenvalues of $A$, then $\ell(A)=2\log r$. Since $r-r^{-1}\le|\tr A|\le r+r^{-1}$, we obtain
\begin{equation}\label{eq:trace-length}
0\le\log(|\tr A|+2)-\frac12\ell(A)\le\log4.
\end{equation}
Consequently, logarithmic traces and translation lengths determine the same limits under any normalization tending to infinity. This gives the description of the Morgan--Shalen boundary in terms of projective length functions \cite[Section~2]{DDW00}.

\begin{definition}\label{def:tree}
A \emph{real tree} is a metric space in which any two points are joined by a unique arc isometric to a real interval. For an isometric action of $\Gamma$ on a real tree $T$, its \emph{marked length function} is defined by
\[
\ell_T(\gamma)=\inf_{x\in T}d_T(x,\gamma x),\quad \gamma\in\Gamma.
\]
The action is \emph{minimal} if it has no proper nonempty invariant subtree. A nontrivial minimal action is \emph{irreducible} if it preserves neither an end nor a line. Its length function is \emph{abelian} if $\ell_T=|c|$ for some homomorphism $c:\Gamma\to\R$.
\end{definition}


A nonzero $q\in Q(X)$ determines a singular flat metric with line element $|q|^{1/2}$. In a natural coordinate $\zeta$ satisfying $q=d\zeta^2$, the horizontal leaves are given by constant values of $\operatorname{Im}\zeta$, with transverse measure $|d\operatorname{Im}\zeta|$. The metric quotient of the foliation lifted to $\widetilde X$ is its dual tree $T_q$, with projection $\pi_q:\widetilde X\to T_q$. We use the standard construction of the dual tree and the continuity of its marked lengths as $q$ varies \cite{HM79,Wol96,FW01}. A flat geodesic is called \emph{nonsingular} if it avoids the zeros of $q$. Its transverse height is the integral of $|\operatorname{Im}\sqrt q|$ along it.

\begin{definition}[Dumas]\label{def:straight}
A continuous map $F:T_q\to T$ is \emph{straight} if its restriction to $\pi_q(J)$ is an isometric embedding whenever $J\subset\widetilde X$ is a nonsingular $q$-geodesic segment. The metric on $\pi_q(J)$ is its transverse metric.
\end{definition}

The following proposition collects consequences of Dumas's results in the normalization used here.

\begin{proposition}\label{prop:dumas-limits}
Suppose that $t_j\to\infty$ and $\varphi_j/t_j^2\to q\ne0$. Then, for every $\gamma\in\Gamma$,
\begin{equation}\label{eq:upper-length}
\ell_{\Hol_X(\varphi_j)}(\gamma)=O_\gamma(t_j).
\end{equation}
Here the implied constant may depend on $\gamma$. If $\delta$ is represented by a nonsingular periodic $q$-geodesic with transverse height $h_q(\delta)$, then
\begin{equation}\label{eq:periodic}
\frac{\ell_{\Hol_X(\varphi_j)}(\delta)}{t_j}\longrightarrow\sqrt2\,h_q(\delta).
\end{equation}
Moreover, every nonzero limit of $\ell_{\Hol_X(\varphi_j)}/t_j$ is represented by a tree admitting a surjective $\Gamma$-equivariant straight map from $T_{2q}$. If all zeros of $q$ are simple, the projective length limit is $[\ell_{T_q}]$.
\end{proposition}

\begin{proof}
The upper bound follows from \cite[Theorem~5.3]{Dum17}, since $\|\varphi_j\|_1^{1/2}/t_j\to\|q\|_1^{1/2}$. For the periodic asymptotic, put $\widehat q=q/\|q\|_1$. By \cite[Theorem~4.3, (4.5)]{Dum17},
\[
\frac{\ell_{\Hol_X(\varphi_j)}(\delta)}{\sqrt{2\|\varphi_j\|_1}}\longrightarrow h_{\widehat q}(\delta)=\frac{h_q(\delta)}{\sqrt{\|q\|_1}},
\]
which gives \eqref{eq:periodic}. The existence of the straight map follows from \cite[Theorem~6.2, Lemma~6.3 and Section~6.6]{Dum17}. Rescaling the source and target metrics by $\sqrt{2\|q\|_1}$ gives a map from $T_{2q}$ to a tree whose length function is the prescribed limit. The assertion for simple zeros follows from \cite[Theorem~B]{Dum17}.
\end{proof}

\begin{proposition}\label{prop:scale-comparison}
Suppose that $\|\varphi_j\|_1\to\infty$ and $\varphi_j/\|\varphi_j\|_1\to q$ with $\|q\|_1=1$. Put $t_j=\sqrt{\|\varphi_j\|_1}$ and $\ell_j=\ell_{\Hol_X(\varphi_j)}$. If $a_j>0$ and $a_j\ell_j\to\ell\ne0$ pointwise on $\Gamma$, then for every nonsingular periodic $q$-geodesic $\delta$ of positive transverse height $h_q(\delta)$,
\begin{equation}\label{eq:scale-comparison}
a_jt_j\longrightarrow\frac{\ell(\delta)}{\sqrt2\,h_q(\delta)}>0.
\end{equation}
Consequently, $\ell_j/t_j$ converges pointwise on $\Gamma$ to a nonzero length function representing the same projective boundary point as $\ell$.
\end{proposition}

\begin{proof}
By Masur's density theorem \cite[Theorem~2]{Mas86}, a nonsingular periodic $q$-geodesic $\delta$ with positive transverse height exists. If $\ell(\delta)=0$, then \eqref{eq:periodic} gives $a_jt_j\to0$. It follows from \eqref{eq:upper-length} that $a_j\ell_j(\gamma)\to0$ for every $\gamma\in\Gamma$, contradicting $\ell\ne0$. Thus $\ell(\delta)>0$, and \eqref{eq:periodic} gives \eqref{eq:scale-comparison}. Since $a_j\ell_j\to\ell$ and $a_jt_j$ converges to a positive constant, the normalized lengths $\ell_j/t_j$ converge pointwise to a positive multiple of $\ell$.
\end{proof}

The matrix estimates needed below are stronger than the length bound and follow directly from the scalar equation.

\begin{lemma}\label{lem:balanced}
Suppose that $t_j\ge1$ and $\varphi_j/t_j^2$ remains in a fixed bounded subset of $Q(X)$. Fix finitely many coordinate disks, compact subsets of these disks, and paths connecting them. In balanced jet coordinates $(y,y'/t_j)$, the fundamental matrices of \eqref{eq:scalar}, their inverses, and their holomorphic derivatives of any fixed order have norms bounded by $e^{Ct_j}$ for suitable constants. The same bounds hold in fixed jet frames.
\end{lemma}

\begin{proof}
Since $Q(X)$ is finite-dimensional, the hypothesis gives a bound of the form $|P_0+\varphi_j|\le Ct_j^2$ on the chosen compact sets. In balanced jet coordinates, the vector $Y=(y,y'/t_j)^T$ satisfies
\begin{equation}\label{eq:balanced}
Y'=\begin{pmatrix}0&t_j\\-(P_0+\varphi_j)/(2t_j)&0\end{pmatrix}Y.
\end{equation}
The coefficient matrix is therefore $O(t_j)$, and Gronwall's inequality gives an $e^{Ct_j}$ bound for transport along each fixed path. The same argument applied to the reversed path gives the bound for inverse transport.

Under a coordinate change $w=w(z)$, the relation $y_w=\sqrt{w'}y_z$ gives the balanced transition matrix
\begin{equation}\label{eq:balanced-transition}
\begin{pmatrix}\sqrt{w'}&0\\ w''/(2t_j(w')^{3/2})&1/\sqrt{w'}\end{pmatrix}.
\end{equation}
Both this matrix and its inverse are uniformly bounded on fixed compact overlaps. Since only finitely many coordinate changes and paths are involved, the preceding transport estimates give the required exponential bounds throughout the chosen collection.

To estimate holomorphic derivatives, we apply Cauchy's estimates on slightly larger disks, where the same exponential bounds hold. Finally, passing from balanced jet coordinates to fixed jet frames introduces only powers of $t_j$, which can be absorbed into the exponential bounds.
\end{proof}


\subsection{Directional fibers and harmonic maps}\label{subsec:fibers}

We use $\mathcal L_X(q)$ as defined in Definition~\ref{def:intro-fiber}, with the $L^1$ norm $\|\cdot\|_1$ on $Q(X)$. The next proposition uses the properness of the holonomy map and Proposition~\ref{prop:dumas-limits} to establish some basic topological properties of these fibers.

\begin{proposition}\label{prop:fiber-topology}
For every $q\in Q(X)$ with $\|q\|_1=1$, the fiber $\mathcal L_X(q)$ is nonempty, compact, and connected, and contains $[\ell_{T_q}]$. The graph of the map $q\mapsto\mathcal L_X(q)$ is closed over the unit sphere in $Q(X)$. Moreover, every fiber that is not a singleton has cardinality $2^{\aleph_0}$.
\end{proposition}

\begin{proof}
The unit sphere of $Q(X)$ with respect to $\|\cdot\|_1$ is radially homeomorphic to a Euclidean sphere. We may therefore choose a nested sequence of connected neighborhoods $V_n$ of $q$ in the unit sphere, with diameters tending to zero. Consider the closures of
\[
\Hol_X\bigl(\{rv:r>n,\ v\in V_n\}\bigr)
\]
in the Morgan--Shalen compactification. Since $V_n$ is connected, these closures are nonempty, compact, and connected. They are also nested, so their intersection is nonempty, compact, and connected. A diagonal argument shows that this intersection is precisely the directional cluster set at $q$. Moreover, the properness of $\Hol_X$ \cite[Theorem~11.4.1]{GKM00}, in the effective form given by \cite[Theorem~5.2]{Dum17}, ensures that the intersection lies entirely in the Morgan--Shalen boundary. Hence $\mathcal L_X(q)$ is nonempty, compact, and connected. Taking the closure of the graph of $\Hol_X$ in the product of the radial compactification of $Q(X)$ and the Morgan--Shalen compactification also shows that the graph of $q\mapsto\mathcal L_X(q)$ is closed.

To show that $[\ell_{T_q}]\in\mathcal L_X(q)$, choose $q_n\to q$ in $Q(X)$ such that $\|q_n\|_1=1$ and all zeros of $q_n$ are simple. By the properness of $\Hol_X$ and Proposition~\ref{prop:dumas-limits}, $\Hol_X(rq_n)$ converges projectively to $[\ell_{T_{q_n}}]$ as $r\to\infty$. The measured foliations and their marked lengths depend continuously on the quadratic differential \cite{HM79,Wol96}. We may therefore choose $r_n\to\infty$ such that $\Hol_X(r_nq_n)\to[\ell_{T_q}]$, proving the desired inclusion.

Finally, suppose that $\mathcal L_X(q)$ contains two distinct points. The distance from one of these points defines a continuous function whose image contains a nondegenerate interval, since $\mathcal L_X(q)$ is connected. Thus the fiber has at least continuum many points. On the other hand, a marked length function is determined by its values on the countable group $\Gamma$, so there are at most continuum many projective marked length functions. It follows that $|\mathcal L_X(q)|=2^{\aleph_0}$.
\end{proof}

We next show how a straight map gives rise to a harmonic map, allowing us to apply rigidity results for marked length functions. We use Korevaar--Schoen's results on the Dirichlet problem and convexity \cite{KS93}, together with the results on the factorization of harmonic maps and minimality in \cite[Proposition~2.4 and the following remark]{DDW00}. 
For a harmonic map $u$ from $\widetilde X$ to a real tree, we write $\Hopf(u)$ for its Hopf differential, defined as the $(2,0)$ part of the pullback metric \cite{KS93}.

\begin{lemma}\label{lem:straight-harmonic}
Let $q\in Q(X)$ be nonzero, and let $F:T_q\to T$ be a continuous, $\Gamma$-equivariant straight map onto a real tree $T$ whose length function is nonzero. Then $u=F\circ\pi_q$ is harmonic, its image is a minimal subtree of $T$, and
\begin{equation}\label{eq:hopf-convention}
\Hopf(u)=-\frac14q.
\end{equation}
Consequently, an irreducible projective marked length function determines the positive real ray of $q$ whenever such a straight map exists.
\end{lemma}

\begin{proof}
On a disk with natural coordinate $\zeta$ that contains no zeros of $q$, straightness implies that $u$ is the composition of $\operatorname{Im}\zeta$ with an isometry onto an interval in $T$. Projecting a competing map onto this interval cannot increase its energy. Since $\operatorname{Im}\zeta$ is harmonic, so is $u$. Direct differentiation gives $\Hopf(\operatorname{Im}\zeta)=-d\zeta^2/4$.

Now let $D$ be a small disk centered at a zero of $q$, with no zeros on its boundary. The map $u$ is continuous on $D$, its energy density is bounded by a constant multiple of $|q|$, and its boundary values are locally Lipschitz. We pass to the metric completion of $T$. The closed convex hull of $u(\partial D)$ is compact. Indeed, the convex hulls of finite nets in $u(\partial D)$ are finite trees that approximate the full hull in the Hausdorff metric. By the Dirichlet theorem for complete nonpositively curved targets \cite[Theorem~2.2 and Section~2.7]{KS93}, there is a harmonic map $v$ with the same boundary values as $u$. On the punctured disk, the function $d_T(u,v)^2$ is bounded and subharmonic \cite[Remark~2.4.3]{KS93}. Its continuous extension across the zero remains subharmonic and vanishes on the boundary. The maximum principle therefore gives $u=v$. Hence $u$ is harmonic across the zero, and \eqref{eq:hopf-convention} holds throughout $X$.

By convexity of the energy, $u$ minimizes energy among equivariant maps on the compact quotient. Its image is minimal by the remark following Proposition~2.4 of \cite{DDW00}. For an irreducible action, the marked length function determines the minimal tree up to equivariant isometry and rescaling \cite[Theorem~3.7]{CM87}, while the normalized Hopf differential on $X$ is unique \cite[Propositions~3.1 and~3.5]{DDW00}. Rescaling the metric on the target by a positive factor multiplies the Hopf differential by the square of that factor. Together with \eqref{eq:hopf-convention}, this shows that the positive real ray of $q$ is determined by the projective marked length function.
\end{proof}


\section{Exclusion of abelian limits}\label{sec:exclusion}

We use the scalar equation \eqref{eq:scalar} and the jet bundle $E$ defined in \eqref{eq:jet}. For a sequence $\varphi_j$ with $\|\varphi_j\|_1\to\infty$, set $t_j=\sqrt{\|\varphi_j\|_1}$ and let $\rho_j$ be a holonomy representation corresponding to $\Hol_X(\varphi_j)$, acting on the fiber $E_{x_0}$ through the associated flat connection.

In this section, we will prove the main theorem that no nonzero abelian length function can occur as a holonomy limit (Theorem~\ref{thm:intro-exclusion}). The proof combines an algebraic approximation with the fact that $M$ has negative degree, which provides the classical degree obstruction in oper theory \cite{Wen16}. Character corrections and elliptic estimates on the oper bundle also appear in \cite[Proposition~6.9]{OSW$^+$24}. Here, estimates for commutators show that the error in equivariance decays faster than any exponential rate and therefore remains negligible under transport along any fixed path.


\subsection{An approximately invariant line}\label{subsec:character-line}

\begin{definition}\label{def:exponential}
Let $t_j\to\infty$. A sequence of vectors or matrices is \emph{exponentially bounded} if its norm is at most $e^{Ct_j}$ for some $C>0$ and all sufficiently large $j$. It is \emph{negligible} if its norm is $O(e^{-At_j})$ for every $A>0$. The same definitions apply to sequences of functions on a compact set, using the specified uniform or Sobolev norm.
\end{definition}

\begin{lemma}\label{lem:character-line}
Suppose that $\rho_j:\Gamma\to\SL_2(\C)$ is a sequence of representations such that $\rho_j(s)$ and $\rho_j(s)^{-1}$ are exponentially bounded for every $s$ in a generating set of $\Gamma$. Assume there exists $\gamma_0\in\Gamma$ such that
\[
\frac{\ell_{\rho_j}(\gamma_0)}{t_j}\longrightarrow2a>0,
\]
and that, for every $\beta\in\Gamma$ and every integer $n\ge1$,
\[
\frac{\ell_{\rho_j}([\gamma_0^n,\beta])}{t_j}\longrightarrow0.
\]
Then, after passing to a subsequence, there exist unit vectors $e_j\in\C^2$ and characters $\chi_j:\Gamma\to\C^*$ such that
\begin{equation}\label{eq:approx-character}
\rho_j(\beta)e_j-\chi_j(\beta)e_j
\end{equation}
is negligible for every $\beta\in\Gamma$. Moreover, the values of $\chi_j$ and their reciprocals are exponentially bounded on every finite set of words.
\end{lemma}

\begin{proof}
Set $A_j=\rho_j(\gamma_0)$ and label its eigenvalues $\lambda_j,\lambda_j^{-1}$ so that $\log|\lambda_j|/t_j\to a$. These eigenvalues are distinct for sufficiently large $j$. The spectral projection
\[
P_j^+=\frac{A_j-\lambda_j^{-1}I}{\lambda_j-\lambda_j^{-1}}
\]
has exponentially bounded norm. Choose unit vectors spanning the two eigenlines and let $S_j$ be the matrix with these vectors as columns. Then $\|S_j\|\le\sqrt2$. The bound on $P_j^+$ also controls the angle between the columns and gives $\|S_j^{-1}\|\le\sqrt2\|P_j^+\|$. Thus the change to this eigenbasis has at most exponential cost.

For $\beta\in\Gamma$, write the matrix of $\rho_j(\beta)$ in this basis as $B_j=\left(\begin{smallmatrix}u_j&b_j\\d_j&v_j\end{smallmatrix}\right)$. Since $\det B_j=1$, we have $u_jv_j-b_jd_j=1$. Direct calculation gives
\begin{equation}\label{eq:commutator}
\tr[A_j^n,B_j]=2-(\lambda_j^n-\lambda_j^{-n})^2b_jd_j.
\end{equation}
By \eqref{eq:trace-length} and the hypothesis, $\log(|\tr[A_j^n,B_j]|+2)=o(t_j)$ for each $n$. It follows that
\[
\limsup_{j\to\infty}\frac{\log|b_jd_j|}{t_j}\le-2na.
\]
Given $A>0$, choose an integer $n$ with $2na>A$. The preceding estimate then shows that $b_jd_j$ is negligible. Consequently, $u_jv_j=1+(\text{negligible})$. Since both diagonal entries are exponentially bounded, they are eventually nonzero and have exponentially bounded reciprocals.

Choose standard handle generators $s_1,\ldots,s_{2g}$ of $\Gamma$. After passing to a subsequence, we may assume that every upper-right and lower-left entry of their matrices has a logarithmic rate in $[-\infty,C]$, where $\log0=-\infty$. For each generator, at least one of these rates is $-\infty$, so the corresponding entry is negligible. Moreover, every entry with a finite rate has an exponential lower bound along this subsequence.

Suppose that a generator matrix $B$ has a nonnegligible upper-right entry, while another generator matrix $C$ has a nonnegligible lower-left entry. Their other off-diagonal entries are negligible, and hence
\[
(BC)_{12}=b_Bv_C+(\text{negligible}),\quad (BC)_{21}=v_Bd_C+(\text{negligible}).
\]
Both leading terms have exponential lower bounds, since the nonnegligible entries $b_B,d_C$ and the diagonal entries $v_B,v_C$ do. Neither term can therefore be canceled by a negligible remainder. Thus the product of the two off-diagonal entries of $BC$ is not negligible, contradicting \eqref{eq:commutator} for the word represented by $BC$. It follows that all nonnegligible off-diagonal entries of the generators occur in the same position. Interchanging the eigenvectors if necessary, we may assume that every lower-left entry is negligible.

Define $\chi_j(s_i)$ to be the upper-left entry of $\rho_j(s_i)$ in the chosen basis. These assignments extend to a character $\chi_j:\Gamma\to\C^*$ because the defining surface relation is a product of commutators. The approximation also holds for inverses. Indeed, for $B=\left(\begin{smallmatrix}u&b\\d&v\end{smallmatrix}\right)$,
\[
B^{-1}\binom10-u^{-1}\binom10=\binom{v-u^{-1}}{-d}=\binom{bd/u}{-d}.
\]
Both entries on the right are negligible. Multiplication then extends the approximation to every word. Returning to the original basis introduces only exponentially bounded factors, which preserve negligibility. Taking $e_j$ to be the corresponding unit eigenvector proves \eqref{eq:approx-character}, together with the asserted bounds on the values of $\chi_j$ and their inverses. This finished the proof. 
\end{proof}


\subsection{Proof of Theorem~\ref{thm:intro-exclusion}}

We now work with holomorphic solutions of the scalar equation \eqref{eq:scalar}. Choose smooth metrics on $X$ and $M$, and write $\|\cdot\|_2$ for the resulting $L^2$ norm. When a unitary line bundle factor is present, we use its flat unitary metric.

\begin{lemma}\label{lem:unitary-gap}
There exists a constant $C_X>0$ such that, for every unitary flat line bundle $U$ on $X$ and every smooth section $s$ of $M\otimes U$,
\begin{equation}\label{eq:unitary-gap}
\|s\|_2\le C_X\|\dbar_{M\otimes U}s\|_2.
\end{equation}
\end{lemma}

\begin{proof}
Since $\deg(M\otimes U)=1-g<0$, the Dolbeault operator on $M\otimes U$ has trivial kernel. The unitary characters of $\Gamma$ form a compact torus. Their flat connections can be written as $d+i\beta$, where $\beta$ is a real harmonic form considered modulo the lattice of integral periods. We may choose representatives in a compact fundamental parallelepiped. On the underlying smooth line bundle of $M$, the corresponding Dolbeault operators take the form $D_\beta=\dbar_M+i\beta^{0,1}$ and satisfy uniform elliptic estimates.

Suppose that \eqref{eq:unitary-gap} fails. By compactness, there exist $\beta_j\to\beta_\infty$ and sections $s_j$ such that $\|s_j\|_2=1$ and $\|D_{\beta_j}s_j\|_2\to0$. The elliptic estimates give a uniform bound for $s_j$ in the first-order Sobolev space $H^1$. By Rellich compactness, a subsequence converges strongly in $L^2$ to a section $s$ with $\|s\|_2=1$. Since $\beta_j\to\beta_\infty$ and $\|D_{\beta_j}s_j\|_2\to0$, we obtain $D_{\beta_\infty}s=0$ weakly. Elliptic regularity then shows that $s$ is a nonzero holomorphic section of a line bundle of negative degree, a contradiction.
\end{proof}

\begin{proposition}\label{prop:oper-line-obstruction}
Suppose that $\varphi_j/t_j^2$ is bounded in $Q(X)$, and let $\rho_j$ denote the corresponding jet monodromy representations. There do not exist unit vectors $e_j$ in the balanced fiber at $x_0$ and exact characters $\chi_j:\Gamma\to\C^*$ whose values and inverses are exponentially bounded on the handle generators, such that \eqref{eq:approx-character} is negligible for every $\beta\in\Gamma$.
\end{proposition}

\begin{proof}
Continue the flat jet with initial value $e_j$ on the universal cover and project it to $M$, obtaining a scalar solution $y_j$. By Lemma~\ref{lem:balanced}, the approximate equivariance in \eqref{eq:approx-character} gives an error that is uniformly negligible on compact subsets of any chosen finite collection of lifted coordinate disks. In particular, continuation around each overlap word differs from multiplication by $\chi_j$ by a solution whose scalar values are negligible on the corresponding compact sets.

The map $\gamma\mapsto\log|\chi_j(\gamma)|$ is a homomorphism to $\R$. Let $\alpha_j$ be its real harmonic representative and set $\eta_j=2\alpha_j^{1,0}$. Then $\eta_j$ is holomorphic and $\operatorname{Re}\eta_j=\alpha_j$. The exponential bounds for $\chi_j$ and its inverses on a chosen homology basis imply
\[
\|\alpha_j\|_{C^k}+\|\eta_j\|_{C^k}\le C_kt_j
\]
for every $k$. On the universal cover, define
\begin{equation}\label{eq:holomorphic-correction}
g_j(x)=\exp\left(\int_{x_0}^x\eta_j\right),\quad
\upsilon_j(\gamma)=\chi_j(\gamma)\exp\left(-\int_\gamma\eta_j\right).
\end{equation}
By construction, $\upsilon_j$ is a unitary character. Moreover, $g_j$ and $g_j^{-1}$ are exponentially bounded on the chosen compact lifts, since these can be reached by paths of bounded length. Consequently, the quotients $y_j/g_j$ are $\upsilon_j$-equivariant up to negligible errors.

Choose a finite good cover $\{V_a\}$ of $X$, lifts of its members, and a smooth partition of unity $\{\psi_a\}$ subordinate to this cover. Let $U_j$ be the flat unitary line bundle that cancels the multiplier $\upsilon_j$, and regard the local quotients as holomorphic sections $f_{a,j}$ of $M\otimes U_j$. Write $T_{ab,j}$ for the transition map of $M\otimes U_j$ from the $b$-frame to the $a$-frame. Then the differences
\[
e_{ab,j}=T_{ab,j}f_{b,j}-f_{a,j}
\]
are negligible on the compact overlap sets used by the partition. Indeed, only finitely many deck words occur, the transition maps for $M$ are bounded on these compact sets, and the unitary transition maps have absolute value one.

Using the partition of unity, define a smooth section $s_j$ by its local expressions
\begin{equation}\label{eq:section-patching}
s_{a,j}=\sum_b\psi_bT_{ab,j}f_{b,j}=f_{a,j}+\sum_b\psi_be_{ab,j}.
\end{equation}
The cocycle identity shows that these expressions agree on overlaps, so $s_j$ is globally defined. Since the local sections are holomorphic and $\sum_b\dbar\psi_b=0$, we obtain
\begin{equation}\label{eq:section-defect}
\dbar s_{a,j}=\sum_b(\dbar\psi_b)e_{ab,j}.
\end{equation}
It follows that $\|\dbar_{M\otimes U_j}s_j\|_2$ is negligible. Furthermore, \eqref{eq:section-patching} shows that $s_j$ differs from each local model by a uniformly negligible error.

We next establish a lower bound for $\|s_j\|_2$. In a coordinate $z$ centered at $x_0$, the normalization of the initial vector in balanced coordinates gives
\[
|y_j(0)|^2+|y_j'(0)/t_j|^2=1.
\]
Choose $r>0$ and put $Y_j(w)=y_j(w/t_j)$. For sufficiently large $j$, the disk $|w|<r$ lies in the rescaled coordinate chart. Since holomorphic Taylor monomials are orthogonal on a disk,
\begin{align}\label{eq:taylor-mass}
\int_{|w|<r}|Y_j(w)|^2\,dA(w)
&\ge\pi r^2|Y_j(0)|^2+\frac{\pi r^4}{2}|Y_j'(0)|^2\\
&\ge\pi\min(r^2,r^4/2).\nonumber
\end{align}
Rescaling the area element gives $\|y_j\|_{L^2(|z|<r/t_j)}\ge c/t_j$. By \eqref{eq:holomorphic-correction}, we have $g_j(0)=1$ and
\[
\left|\int_0^z\eta_j\right|\le Ct_j|z|\le Cr
\]
on this disk. Hence both $g_j$ and its inverse are uniformly bounded there, and dividing by $g_j$ preserves the lower bound up to a constant independent of $j$. The chosen smooth metrics affect the estimate only by such constants. Since the error in \eqref{eq:section-patching} is negligible, we conclude that
\begin{equation}\label{eq:mass-lower}
\|s_j\|_2\ge c'/t_j
\end{equation}
for all sufficiently large $j$.

Applying \eqref{eq:unitary-gap} to $s_j$ now contradicts the negligible bound in \eqref{eq:section-defect} and the lower bound in \eqref{eq:mass-lower}. This proves the proposition.
\end{proof}

\begin{proof}[Proof of Theorem~\ref{thm:intro-exclusion}]
Suppose that $\|\varphi_j\|_1\to\infty$ and that the corresponding projective marked length functions converge to $[|c|]$ for some nonzero homomorphism $c:\Gamma\to\R$. After passing to a subsequence on which the normalized Schwarzians converge, \eqref{eq:scale-comparison} shows that the lengths normalized by $t_j=\sqrt{\|\varphi_j\|_1}$ converge pointwise to a positive multiple of $|c|$. Rescaling $c$ by a positive constant, we may therefore assume that $\ell_{\rho_j}(\gamma)/t_j\to|c(\gamma)|$ for every $\gamma\in\Gamma$.

Choose $\gamma_0\in\Gamma$ with $c(\gamma_0)\ne0$. Since $c$ is a homomorphism, every commutator $[\gamma_0^n,\beta]$ has zero limiting length. Lemma~\ref{lem:balanced} gives the exponential matrix bounds required by Lemma~\ref{lem:character-line}, which then yields an approximate character line. This contradicts Proposition~\ref{prop:oper-line-obstruction} and completes the proof.
\end{proof}


\subsection{Finite covers and commutator bounds}\label{subsec:commutators}

\begin{corollary}\label{cor:irreducible-limits}
Every nonzero length function arising as a holonomy boundary limit on $X$ restricts to a nonzero, nonabelian length function on every finite-index subgroup of $\Gamma$. Its minimal tree action is irreducible.
\end{corollary}

\begin{proof}
Let $H$ be a finite-index subgroup of $\Gamma$, and let $p:Y\to X$ be the corresponding connected unramified cover. Pulling back a projective structure from $X$ to $Y$ restricts its holonomy representation to $H$. Moreover, $g(Y)=1+\deg(p)(g(X)-1)\ge2$ and $K_Y=p^*K_X$. The pullback map $p^*:Q(X)\to Q(Y)$ satisfies $\|p^*\varphi\|_{1,Y}=\deg(p)\|\varphi\|_{1,X}$ and therefore preserves divergence. The restriction of a nonzero length function to $H$ is also nonzero, since some positive power of an element with positive length belongs to $H$. Applying Theorem~\ref{thm:intro-exclusion} to $Y$ therefore shows that the restricted length function cannot be abelian.

We next consider the minimal tree action. If it fixes an end, its length function is the absolute value of its Busemann character \cite[Section~2]{CM87}. If it preserves a line, the subgroup preserving both ends has index at most two and admits a nonzero translation character. Both cases contradict the first assertion. This proves the assertion.
\end{proof}

\begin{corollary}\label{cor:finite-commutators}
There exist finitely many commutators $w_1,\ldots,w_k\in\Gamma$ and constants $\varepsilon,R>0$ such that
\[
\max_{1\le i\le k}\ell_{\Hol_X(\varphi)}(w_i)\ge\varepsilon\sqrt{\|\varphi\|_1}\quad\text{whenever }\|\varphi\|_1>R.
\]
\end{corollary}

\begin{proof}
Enumerate the commutators of $\Gamma$ and suppose the conclusion is false. For each $j$, choose $\varphi_j\in Q(X)$ with $\|\varphi_j\|_1>j$ such that the lengths of the first $j$ commutators are all less than $t_j/j$, where $t_j=\sqrt{\|\varphi_j\|_1}$. Passing to a subsequence, we may assume that $\varphi_j/\|\varphi_j\|_1\to q$ for some $q\in Q(X)$ with $\|q\|_1=1$. The bounds for each element of $\Gamma$ allow us to extract a further subsequence along which $\ell_{\Hol_X(\varphi_j)}(\gamma)/t_j$ converges for every $\gamma\in\Gamma$. By \eqref{eq:periodic}, at least one of these limits is positive, whereas every commutator has limiting length zero. The exponential bounds from Lemma~\ref{lem:balanced} allow us to apply Lemma~\ref{lem:character-line} to an element of positive limiting length. The resulting approximate equivariance contradicts Proposition~\ref{prop:oper-line-obstruction}.
\end{proof}


\section{Nonuniqueness at multiple zeros}\label{sec:nonuniqueness}

Choose a nonzero $q\in Q(X)$ with a multiple zero and a reference projective connection $P_0$. In this section, we prove Theorem~\ref{thm:intro-nonunique} by constructing globally holomorphic perturbations of $t^2q$ and prescribing the exact holonomy trace of a chosen element of $\Gamma$.


\subsection{Closed geodesics and polynomial equations}\label{subsec:closed-geodesics}

We use the flat metric $|q|$ on $X$. A zero of $q$ of order $n$ has cone angle $(n+2)\pi$. The cone-angle criterion then shows that the universal cover of $X$, with the lifted path metric, is a proper, complete, and geodesically complete $\operatorname{CAT}(0)$ space; see \cite[Chapter~II]{BH99} and \cite[Section~2]{Dum17}.

\begin{lemma}\label{lem:prescribed-turn}
Let $p$ be a zero of $q$, and let $J_-$ and $J_+$ be open intervals in the circle of directions at $p$. Suppose that any two directions, one from each interval, divide the circle into two arcs, each of length strictly greater than $\pi$. Then there exists a closed $q$-geodesic passing through $p$ with one outward direction in $J_-$ and the other in $J_+$.
\end{lemma}

\begin{proof}
Let $Y$ be the universal cover of $X$ with the lifted flat metric, and choose a lift $\widetilde p$ of $p$. The cocompact action of $\Gamma$ makes $Y$ quasi-isometric to a Cayley graph of $\Gamma$, so $Y$ is Gromov hyperbolic. Choose two geodesic rays from $\widetilde p$ whose initial directions lie in $J_-$ and $J_+$, respectively. Geodesic completeness allows them to continue across any further cone points. Since both angles between their initial directions exceed $\pi$, their union is a local geodesic and hence a complete geodesic line. In particular, its two points at infinity are distinct.

Let $\partial Y$ denote the visual boundary of $Y$, consisting of asymptotic classes of geodesic rays. The initial direction of a geodesic ray from $\widetilde p$ depends continuously on its endpoint in $\partial Y$. Indeed, convergence in the visual boundary is equivalent to compact-open convergence of geodesic rays based at $\widetilde p$. Their intersections with a sufficiently small circle centered at $\widetilde p$, enclosing no other singularities, therefore converge, as do their initial directions. Consequently, there exist disjoint nonempty open sets $V_-,V_+\subset\partial Y$ such that the rays from $\widetilde p$ toward points of $V_-$ and $V_+$ have initial directions in $J_-$ and $J_+$, respectively.

The pairs of repelling and attracting fixed points of nontrivial elements of $\Gamma$ are dense in the complement of the diagonal in $\partial Y\times\partial Y$; see \cite[Lemma~2.2]{BM12}. We may therefore choose $\gamma\in\Gamma$ whose repelling and attracting fixed points belong to $V_-$ and $V_+$, respectively. The geodesic rays from $\widetilde p$ toward these points join to form a geodesic line $a$. By the flat-strip theorem \cite[Theorem~II.2.13]{BH99}, any other geodesic line with the same points at infinity would bound a flat strip together with $a$. Such a strip cannot have positive width, since it would have angle $\pi$ on one side of $a$ at $\widetilde p$, whereas both angles there are strictly greater than $\pi$. Thus $a$ is the unique geodesic line joining these points at infinity. Since $\gamma$ preserves both points, it also preserves $a$. Moreover, $\gamma$ acts on $a$ by a nonzero translation, because a nonidentity deck transformation has no fixed point in $Y$. The projection of $a$ to $X$ is therefore the required closed $q$-geodesic.
\end{proof}

\begin{definition}
In a local coordinate satisfying $q=z^n\,dz^2$, the $n+2$ open regions between consecutive horizontal rays of $q$ are called \emph{Stokes sectors}. Each sector has flat angle $\pi$. Two sectors are \emph{adjacent} if their closures share a boundary ray. For a scalar polynomial equation, a \emph{recessive solution} in a sector is a nonzero solution that decays exponentially along rays in the interior of that sector.
\end{definition}

Let $m\ge2$ be the maximum order of the zeros of $q$, and choose a zero $p$ of order $m$. The central directions of two nonadjacent Stokes sectors at $p$ divide the circle of directions into two arcs, each of flat length at least $2\pi$. Sufficiently small intervals around these directions therefore satisfy the hypotheses of Lemma~\ref{lem:prescribed-turn}, which gives a closed $q$-geodesic $\gamma$. A period of its axis contains only finitely many zeros, since the lifted zero set is locally finite and the period is compact. The same zero may be visited more than once.

We first assume that no saddle edge of $\gamma$ is horizontal. At each zero visited by $\gamma$, the two outward directions must lie in distinct Stokes sectors. Otherwise, they would enclose a flat angle less than $\pi$, contradicting the local geodesic criterion. At $p$, these directions lie in nonadjacent sectors by construction. We remove the assumption on horizontal saddle edges in Section~\ref{subsec:trace-prescription}.

We choose $\eta\in Q(X)$ whose local coefficient at $p$ is $1$ in a coordinate satisfying $q=z^m\,dz^2$. Such a differential exists because $K_X^2$ is generated by its global sections. Indeed, Riemann--Roch and $H^0(X,K_X^{-1}(p))=0$ give $h^0(K_X^2)-h^0(K_X^2(-p))=1$. We define
\begin{equation}\label{eq:global-perturbation}
b=\frac4{m+2},\quad \varphi_t(a)=t^2q+t^b a\eta,\quad t>0,\quad a\in\C.
\end{equation}
At a zero of order $n\le m$ visited by $\gamma$, we choose a local coordinate in which $q=z^n\,dz^2$ and write $\eta=(\sum_{k\ge0}c_kz^k)\,dz^2$. Under the rescaling $z=\epsilon_n w$, where $\epsilon_n=t^{-2/(n+2)}$, the scalar equation has leading coefficient $w^n/2$. The coefficient of $a c_kw^k/2$ is
\begin{equation}\label{eq:inner-exponent}
t^{b-2(k+2)/(n+2)},\quad b-\frac{2(k+2)}{n+2}\le0.
\end{equation}
The exponent vanishes precisely when $n=m$ and $k=0$. The reference coefficient contributes $\epsilon_n^2P_0(\epsilon_n w)/2$. Consequently, the local equations converge uniformly on compact subsets of the $(w,a)$-space to
\begin{equation}\label{eq:inner-model}
U''+\frac12(w^n+a c)U=0,\quad c=\begin{cases}c_0,&n=m,\\0,&n<m.\end{cases}
\end{equation}
The maximality of $m$ ensures a bounded polynomial model at every zero visited by $\gamma$.

We need two classical facts about the recessive solutions of \eqref{eq:inner-model}. The first concerns the selected zero, where $n=m$ and $c=1$. For certain isolated values of $a$, the recessive lines associated with two nonadjacent sectors coincide. When $m\ge3$, this follows from Shin's theorem on polynomial eigenvalue problems \cite[Theorem~1.1 and Section~2]{Shi04}. Indeed, the substitution $w=c_*x$, where $c_*^{m+2}=-2$, transforms the equation into
\begin{equation}\label{eq:schrodinger-model}
-u''+x^m u=\lambda u,\quad \lambda=\frac{a c_*^2}{2}.
\end{equation}
To apply Shin's theorem, let $d$ denote the cyclic separation between the two sectors. Since they are nonadjacent, we have $2\le d\le m$, and hence $\ell=d-1$ lies in the range $1\le\ell\le m-1$ required by the theorem. The rotations described in \cite{Shi04} transform the boundary rays into those of the sectors indexed by $0$ and $-d$ for the monic equation. Further rotations by roots of unity cover every cyclic pair. The spectrum is discrete and infinite. When $m=2$, the same conclusion follows from the classical Weber equation. Indeed, $e^{-x^2/2}$ solves \eqref{eq:schrodinger-model} with $\lambda=1$ and is recessive along both real rays, while rotation by $i$ handles the other nonadjacent pair \cite[Chapter~12]{OOL$^+$26}. An isolated zero of finite positive multiplicity suffices for the argument that follows.

Second, when $a=0$, the recessive lines associated with any two distinct sectors are linearly independent. This follows from the classical representation in terms of Bessel functions \cite[Chapter~10]{OOL$^+$26}. For the equation $-u''+x^n u=0$, set $\nu=1/(n+2)$ and $\rho=(n+2)/2$. A solution that is recessive along the positive real ray is given by
\[
f(x)=\sqrt{x}\,K_\nu(x^\rho/\rho)=A+B x+O(x^{n+2}),
\]
where
\[
A=\frac{1}{2}\Gamma(\nu)(n+2)^\nu,\quad B=\frac{1}{2}\Gamma(-\nu)(n+2)^{-\nu}.
\]
Both $A$ and $B$ are nonzero. Now let $\omega=e^{2\pi i/(n+2)}$. The rotated solutions $f_k(x)=f(\omega^{-k}x)$ have initial slopes $\omega^{-k}B/A$, which are pairwise distinct. Hence they span distinct lines in the two-dimensional solution space. A constant complex rescaling gives the same conclusion for \eqref{eq:inner-model} when $a=0$.


\subsection{Uniform estimates for connection matrices near a resonance}\label{subsec:local-matrices}

We next establish an estimate that allows us to incorporate the local spectral zero into the calculation of the global holonomy trace. We consider an equation of the form \eqref{eq:inner-model} and two rays in the interiors of its chosen sectors. Before deriving the estimates, we choose a closed parameter disk $D\subset\C$ containing both $a=0$ and a selected spectral value. The constructions below depend holomorphically on $a$ in a neighborhood of $D$.

Along an outward ray $w=e^{i\theta}r$, we write the equation as $U_{rr}=Q(r,a)U$ and choose a square root $p$ of $Q$ such that
\begin{equation}\label{eq:radial-root}
\Re p(r,a)\ge c_\theta r^{n/2}>0\quad(r\ge R_0,\ a\in D).
\end{equation}
By choosing $R_0$ sufficiently large, we can ensure that this inequality holds uniformly for $a\in D$. We then set 
\[
\Phi(r,a)=\int_{R_0}^r p(s,a)\,ds, \quad \nu_n=\frac{n+2}{2}.
\]

\begin{lemma}\label{lem:sector-bases}
On each chosen ray there is a holomorphic family of recessive solutions $f$ and complementary solutions $g$ with constant nonzero Wronskian such that, uniformly for $a\in D$,
\begin{align}\label{eq:sector-asymptotics}
f&=p^{-1/2}e^{-\Phi}(1+O(r^{-\nu_n})),& f_r&=-p^{1/2}e^{-\Phi}(1+O(r^{-\nu_n})),\\
g&=p^{-1/2}e^{\Phi}(1+O(r^{-\nu_n})),& g_r&=p^{1/2}e^{\Phi}(1+O(r^{-\nu_n})).\notag
\end{align}
The solutions extend to the whole local complex plane and depend holomorphically on $a$.
\end{lemma}

\begin{proof}
The existence of sectorial solutions is classical \cite{Sib75}. We establish the uniform estimates and normalization needed for the subsequent matrix comparison. Using $\xi=\Phi(r,a)$ as the new independent variable and writing $U=p^{-1/2}u$, we obtain
\begin{equation}\label{eq:liouville-residual}
u_{\xi\xi}=(1-\mathfrak g)u,\quad \text{where}\quad \mathfrak g=\frac{3p_r^2}{4p^4}-\frac{p_{rr}}{2p^3}=\frac{5Q_r^2}{16Q^3}-\frac{Q_{rr}}{4Q^2}.
\end{equation}
The leading polynomial term, together with \eqref{eq:radial-root}, yields
\[
\mathfrak g=O(r^{-n-2}),\quad \int_r^\infty |\mathfrak g|\,|d\xi|\le C r^{-\nu_n}.
\]
For $v=e^\xi u$, the recessive normalization leads to the integral equation
\begin{equation}\label{eq:volterra}
v(\xi)=1+\frac12\int_\xi^\infty\bigl(e^{2(\xi-s)}-1\bigr)\mathfrak g(s)v(s)\,ds.
\end{equation}
Along the image of the ray, $\Re\xi$ increases, so the kernel in \eqref{eq:volterra} has modulus at most one. Successive iterated integrals are bounded by powers of the tail integral divided by the corresponding factorials. The resulting series converges uniformly for $a\in D$ and gives $v=1+O(r^{-\nu_n})$. Differentiating \eqref{eq:volterra}, we obtain
\[
v_\xi=\int_\xi^\infty e^{2(\xi-s)}\mathfrak g(s)v(s)\,ds,\quad v_{\xi\xi}-2v_\xi=-\mathfrak g v.
\]
These identities give the required estimates for $f=p^{-1/2}e^{-\Phi}v$ and its derivative. The uniform convergence of the iterates also establishes holomorphic dependence on $a$.

We may increase $R_0$ so that $f$ has no zeros on the outward ray for any $a\in D$. We then define
\begin{equation}\label{eq:complementary-solution}
g(r,a)=2f(r,a)\int_{R_0}^r f(s,a)^{-2}\,ds.
\end{equation}
Writing $W_r(f,g)=f g_r-f_r g$ for the Wronskian, we obtain $W_r(f,g)=2$. Moreover, $f^{-2}=p e^{2\Phi}(1+O(r^{-\nu_n}))$. The integral of the leading term equals $e^{2\Phi}/2$ up to the contribution from the lower limit. To estimate the remaining integral, we split the interval at $r/2$. After division by $e^{2\Phi(r)}$, the contribution from the first part is exponentially small. On the second part, the relative error is $O(r^{-\nu_n})$, and the inequality $|p|\le C\Re p$ bounds its absolute integral by a constant multiple of $e^{2\Re\Phi(r)}$. This proves the asymptotic estimate for $g$, while differentiation of \eqref{eq:complementary-solution} gives the corresponding derivative estimate. Since the differential equation has polynomial coefficients, both solutions extend from their initial data to the entire complex plane. The initial data depend holomorphically on $a$, and so do the extended solutions. This completes the proof.
\end{proof}

For $p\in\C^*$, we define the WKB frame
\begin{equation}\label{eq:wkb-frame}
F(p)=\begin{pmatrix}p^{-1/2}&p^{-1/2}\\p^{1/2}&-p^{1/2}\end{pmatrix}.
\end{equation}
Its determinant is $-2$. Along the incoming ray, differentiation in the direction of travel is the negative of radial differentiation. Consequently, at radius $M$, the ordered basis of exact solutions $(f_{\rm in},g_{\rm in})$ is represented in the travel WKB frame by
\[
(I+E_{\rm in})\diag(e^{-\Phi_{\rm in}},e^{\Phi_{\rm in}}),\quad \|E_{\rm in}\|\le C M^{-\nu_n}.
\]
Along the outgoing ray, the ordered basis $(g_{\rm out},f_{\rm out})$ is represented by
\[
(I+E_{\rm out})\diag(e^{\Phi_{\rm out}},e^{-\Phi_{\rm out}}),\quad \|E_{\rm out}\|\le C M^{-\nu_n}.
\]
Here $\Phi_{\rm in}$ and $\Phi_{\rm out}$ are the corresponding phase integrals evaluated at radius $M$. The constant half-form factors arising from the change between radial and holomorphic coordinates are included in these bases.

Expressing both pairs of solutions in the same holomorphic Cauchy frame, we define the change-of-basis matrix $B(a)$ by
\[
(f_{\rm in},g_{\rm in})=(g_{\rm out},f_{\rm out})B(a).
\]
Since both pairs form bases with nonzero Wronskians throughout $D$, the matrix $B(a)$ is holomorphic and bounded on $D$. Evaluating all Wronskians in the same coordinate, we obtain
\begin{equation}\label{eq:connection-coefficient}
C_v(a):=B_{11}(a)=\frac{W(f_{\rm in},f_{\rm out})}{W(g_{\rm out},f_{\rm out})}.
\end{equation}
Here the subscript $v$ labels the visit to a zero. The denominator remains nonzero even when the two recessive lines coincide, since $g_{\rm out}$ and $f_{\rm out}$ are linearly independent.

\begin{proposition}\label{prop:local-transfer}
Let $L_{v,M}(a)$ be the transfer across this zero, expressed in the travel WKB frames on the incoming and outgoing rays. Define
\[
N_{v,M}(a)=\exp(\Phi_{\rm in}(M,a)+\Phi_{\rm out}(M,a)),\quad P=\diag(1,0).
\]
Then $N_{v,M}$ is holomorphic and nowhere zero, and
\begin{equation}\label{eq:local-transfer-limit}
N_{v,M}^{-1}L_{v,M}\longrightarrow C_v(a)P
\end{equation}
uniformly for $a\in D$ as $M\to\infty$. Moreover, the matrices $N_{v,M}^{-1}L_{v,M}$ are uniformly bounded on $D$.
\end{proposition}

\begin{proof}
We express the transfer using the two bases of exact solutions. After dividing by $N_{v,M}$, we obtain
\begin{equation}\label{eq:full-local-matrix}
\frac{L_{v,M}}{N_{v,M}}=(I+E_{\rm out})
\begin{pmatrix}
B_{11}&B_{12}e^{-2\Phi_{\rm in}}\\
B_{21}e^{-2\Phi_{\rm out}}&B_{22}e^{-2\Phi_{\rm in}-2\Phi_{\rm out}}
\end{pmatrix}
(I+E_{\rm in})^{-1}.
\end{equation}
The real parts of both phase integrals tend to infinity uniformly for $a\in D$. Since $B$ is bounded and $E_{\rm in},E_{\rm out}\to0$, the right-hand side converges uniformly to $B_{11}(a)P=C_v(a)P$. The same formula also shows that the normalized matrices remain uniformly bounded, even when $B_{11}(a)=0$.
\end{proof}

We apply this construction to each of the finitely many visits of $\gamma$ to a zero. The resulting product
\begin{equation}\label{eq:connection-product}
C(a)=\prod_v C_v(a)
\end{equation}
is holomorphic. Since the recessive lines associated with distinct sectors are independent at $a=0$, every factor $C_v(0)$ is nonzero, and hence $C(0)\ne0$. On the other hand, at a spectral value associated with the selected zero, at least one factor vanishes. Thus $C$ has an isolated zero. Although the values of $\eta$ at different zeros may be correlated and $\gamma$ may visit the same zero more than once, these effects can increase the multiplicity of a zero of $C$ but cannot make $C$ vanish identically.


\subsection{Matching along the closed geodesic}\label{subsec:global-matching}

We now compare the transfers arising from the polynomial models with those of the scalar equation associated with the global differential \eqref{eq:global-perturbation}. Consider an oriented nonsingular saddle edge of flat length $L_e$. Along this edge, we choose a natural coordinate $\zeta=\zeta_0+e^{i\theta_e}s$, where $s$ denotes arclength and $\sin\theta_e\ne0$. In the corresponding scalar half-form coordinate, the equation takes the form $v_{ss}=p_s^2v$, where
\begin{equation}\label{eq:outer-potential}
p_s^2=-\frac{e^{2i\theta_e}}2\bigl(t^2+t^ba h(s)+R_\zeta(s)\bigr),\quad h=\eta/q.
\end{equation}
Here $R_\zeta$ denotes the reference projective connection in the natural coordinate. Near a zero of order $n$, set $\alpha_n=2n/(n+2)$ and $\nu_n=(n+2)/2$, and write $d$ for the flat distance to the zero. Since $\zeta$ is a constant multiple of $z^{(n+2)/2}$ on the chosen sector, the coordinate-change rule gives
\begin{equation}\label{eq:outer-singular-bounds}
h=O(d^{-\alpha_n}),\quad h'=O(d^{-\alpha_n-1}),\quad R_\zeta=O(d^{-2}),\quad R_\zeta'=O(d^{-3}).
\end{equation}
Here derivatives are taken along the edge. We truncate the local disk at $|z|=Mt^{-2/(n+2)}$. The flat distance from the zero to this cutoff is
\begin{equation}\label{eq:cutoff-distance}
d_{\rm cut}=\frac2{n+2}\frac{M^{\nu_n}}t.
\end{equation}

\begin{lemma}\label{lem:outer-transfer}
Let $E_{e,t,M}(a)$ be the transfer along the remaining part of the edge in the frames \eqref{eq:wkb-frame}. There is a square root in \eqref{eq:outer-potential}, holomorphic for $a$ near $D$, such that, with $I_{e,t,M}=\int p_s\,ds$,
\begin{equation}\label{eq:outer-transfer-limit}
e^{-I_{e,t,M}}E_{e,t,M}=P+o_M(1)+o_{t,M}(1).
\end{equation}
The estimates are uniform on $D$, with $M$ chosen first and then $t$ sufficiently large. For each fixed sufficiently large $M$,
\begin{equation}\label{eq:outer-exponent}
\frac{\Re I_{e,t,M}(a)}t\longrightarrow \frac{|\sin\theta_e|L_e}{\sqrt2}
\end{equation}
uniformly on $D$.
\end{lemma}

\begin{proof}
Choose $k_{\theta_e}$ with $k_{\theta_e}^2=-e^{2i\theta_e}/2$ and $\Re k_{\theta_e}=|\sin\theta_e|/\sqrt2$. Near a zero of order $n$, the perturbation in \eqref{eq:outer-potential}, relative to $t^2$, is bounded at the cutoff by
\begin{equation}\label{eq:outer-relative-error}
C\bigl(t^{b-4/(n+2)}M^{-n}+M^{-n-2}\bigr).
\end{equation}
Indeed, substituting \eqref{eq:cutoff-distance} into $t^{b-2}|h|+t^{-2}|R_\zeta|$ gives precisely these powers of $t$ and $M$. Since $n\le m$, we have $b-4/(n+2)\le0$. We may therefore choose $M$ sufficiently large and then take $t$ sufficiently large so that the square root can be chosen near $tk_{\theta_e}$ throughout the truncated edge and over the parameter disk $D$. In particular, $\Re p_s\ge c_e t$.

Differentiating the square root and applying \eqref{eq:outer-singular-bounds}, we obtain
\begin{equation}\label{eq:outer-coupling}
\int\left|\frac{p_s'}{p_s}\right|\,ds\le C\sum_{n\text{ at the two ends}}\bigl(t^{b-4/(n+2)}M^{-n}+M^{-n-2}\bigr)+o_t(1).
\end{equation}
The contribution from the first derivative term is bounded by a constant multiple of $t^{b-2}d_{\rm cut}^{-\alpha_n}$, while the reference term contributes at most a constant multiple of $t^{-2}d_{\rm cut}^{-2}$. Away from the zeros, the remaining contribution tends to zero as $t\to\infty$ for each $M$.

In the frame $F(p_s)$, the first-order equation has diagonal part $\diag(p_s,-p_s)$ and off-diagonal part with norm $|p_s'/2p_s|$. We expand its fundamental matrix using ordered Duhamel integrals. In each term, the product of the diagonal propagation factors has modulus at most $\exp\int\Re p_s\,ds$. Thus, after dividing by $e^{I_{e,t,M}}$, the sum of all terms containing at least one off-diagonal factor is bounded by
\[
\exp\left(\int|p_s'/2p_s|\,ds\right)-1.
\]
The term without any off-diagonal factor is $\diag(1,e^{-2I_{e,t,M}})$. Its second diagonal entry tends to zero, so \eqref{eq:outer-coupling} gives \eqref{eq:outer-transfer-limit} for the entire matrix.

To determine the leading exponent, we use the square-root estimate
\[
|p_s-tk_{\theta_e}|\le C\bigl(t^{b-1}|h|+t^{-1}|R_\zeta|\bigr).
\]
Near a simple zero, the integral of $d^{-\alpha_n}$ remains bounded, and $b\le1$. For $n=2$, the integral grows logarithmically, but its multiplier satisfies $t^{b-1}\le1$. For $n>2$, its contribution is $O_M(t^{b-4/(n+2)})$, which remains bounded. The reference term contributes $O_M(1)$. Removing the two cut portions changes $tk_{\theta_e}L_e$ by $O_M(1)$. We therefore obtain
\[
I_{e,t,M}=tk_{\theta_e}L_e+O_M(1+\log t),
\]
uniformly on $D$, which proves \eqref{eq:outer-exponent}.
\end{proof}

We must also check that the local and outer frames use the same ordering, since this determines which matrix entry contributes to the trace product. On an incoming local ray, we use the travel coordinate $x=-t^{2/(n+2)}|z|$, while on an outgoing ray we use $x=t^{2/(n+2)}|z|$. In both cases, $x$ increases in the direction of travel. Thus, writing $r=dx/ds$, we have
\begin{equation}\label{eq:travel-derivative}
r=t|x|^{-n/2}>0.
\end{equation}
For $U_{xx}=p_x^2U$, set $v(s)=r^{-1/2}U(x(s))$. Direct differentiation gives
\begin{equation}\label{eq:half-form-match}
\binom v{v_s}=H\binom U{U_x},\quad H=\begin{pmatrix}r^{-1/2}&0\\-r'/(2r^{3/2})&r^{1/2}\end{pmatrix},\quad p_s^2=r^2p_x^2-\frac{1}{2}S(x;s),
\end{equation}
where $S(x;s)$ is the Schwarzian derivative and $r'$ denotes $dr/ds$. Since $x$ is a constant multiple of $d^{2/(n+2)}$, we obtain
\[
S(x;s)=\frac{1-(2/(n+2))^2}{2d^2}.
\]
At $|x|=M$, this gives
\begin{equation}\label{eq:half-form-errors}
\frac{|S(x;s)|}{|r^2p_x^2|}=O(M^{-n-2}),\quad \frac{|r'|}{r^2|p_x|}=O(M^{-\nu_n}).
\end{equation}
We choose both $p_s$ and $p_x$ to have positive real parts in their travel coordinates, so that $p_s/(rp_x)\to1$. To compare the frames explicitly, set $A=(p_s/(rp_x))^{1/2}$, $D_0=A^{-1}$, and $E_0=r'/(2r^{3/2}\sqrt{p_sp_x})$. Direct matrix multiplication gives
\begin{equation}\label{eq:frame-comparison}
F(p_s)^{-1}HF(p_x)=\frac12\begin{pmatrix}A+D_0-E_0&A-D_0-E_0\\A-D_0+E_0&A+D_0+E_0\end{pmatrix}=I+O(M^{-\nu_n}).
\end{equation}
Indeed, $A-1=O(M^{-n-2})$ and $E_0=O(M^{-\nu_n})$. The positivity of $r$ in \eqref{eq:travel-derivative} ensures that the growing coordinate on the incoming edge matches the first coordinate of the local transfer. We account for the constant complex half-form factors by using a common holomorphic Cauchy frame at each zero. The remaining spin factor around the closed curve is a fixed sign.

For each $M$, the rescaled coefficients of the global equation converge to those of \eqref{eq:inner-model} uniformly on the local disk and for $a\in D$ as $t\to\infty$. The integral equation for a fundamental matrix then gives convergence of its Cauchy data and transfer matrices, uniformly on $D$ and with holomorphic dependence on $a$. Consequently, the frame comparison for the global equation has error $O(M^{-\nu_n})+o_{t,M}(1)$. The constants may depend on $M$, and $t$ is chosen sufficiently large only after the cut radius has been determined.

\begin{proposition}\label{prop:global-trace}
For the closed geodesic $\gamma$ with nonhorizontal saddle edges, let $H=\sum_e|\sin\theta_e|L_e>0$. There are a parameter disk $D$, a holomorphic function $C$ with $C(0)\ne0$ and an isolated zero in $D$, and holomorphic nowhere-zero functions $G_{t,M}$ near $D$, such that
\begin{equation}\label{eq:global-trace}
\frac{\tr\Hol_X(\varphi_t(a))(\gamma)}{G_{t,M}(a)}=C(a)+o_M(1)+o_{t,M}(1)
\end{equation}
uniformly on $D$. Moreover, for fixed sufficiently large $M$,
\begin{equation}\label{eq:global-growth}
\frac{\log|G_{t,M}(a)|}{t}\longrightarrow\frac{H}{\sqrt2}
\end{equation}
uniformly on $D$.
\end{proposition}

\begin{proof}
We compute the holonomy from a nonsingular point of $\gamma$, expressing it as a finite ordered product of outer transfers, local transfers, and half-form changes. Define $G_{t,M}$ as the product of $e^{I_{e,t,M}}$ over the edges, $N_{v,M}$ over the visits to zeros, and the fixed spin factor. This function is holomorphic and nowhere zero. After dividing each transfer by its corresponding factor, the outer matrices converge to $P$, the local matrices converge to $C_vP$, and the changes of frame converge to $I$. All these matrices are uniformly bounded, even when $C_v=0$. Since $P^k=P$ and $\tr P=1$, a telescoping estimate for the finite product gives \eqref{eq:global-trace}, with $C=\prod_v C_v$. The properties of $C$ follow from the discussion after \eqref{eq:connection-product}. For each $M$, the logarithmic moduli of the local factors remain bounded as $t\to\infty$. Summing \eqref{eq:outer-exponent} over the edges then proves \eqref{eq:global-growth}.
\end{proof}


\subsection{Proof of Theorem~\ref{thm:intro-nonunique}}\label{subsec:trace-prescription}

Let $\Omega\subset D$ be a small closed disk around an isolated zero of $C$, chosen so that $C$ does not vanish on $\partial\Omega$, and set $c_*=\min_{\partial\Omega}|C|>0$. We first take $M$ sufficiently large that the $M$-error in \eqref{eq:global-trace} is less than $c_*/4$, and then take $t$ sufficiently large that the remaining error is less than $c_*/4$. Rouch\'e's theorem then shows that the exact holonomy trace has a zero in $\Omega$.

More generally, let $0<\lambda<\sqrt2H$ and set
\begin{equation}\label{eq:trace-target}
T_t=2\cosh(\lambda t/2).
\end{equation}
By \eqref{eq:global-growth},
\[
\frac1t\log\left|\frac{T_t}{G_{t,M}(a)}\right|\longrightarrow\frac\lambda2-\frac H{\sqrt2}<0
\]
uniformly for $a\in\partial\Omega$. Hence $|T_t/G_{t,M}(a)|<c_*/4$ on $\partial\Omega$ for sufficiently large $t$. Applying Rouch\'e's theorem to $(\tr\Hol_X(\varphi_t(a))(\gamma)-T_t)/G_{t,M}(a)$, we obtain $a_t\in\Omega$ satisfying the exact trace equation. The argument also accounts for multiple zeros of $C$, counted with their multiplicities. Since $a_t$ remains bounded, we have $\varphi_t(a_t)/t^2\to q$.

A matrix in $\SL_2(\C)$ with trace zero has eigenvalues $i$ and $-i$, so its translation length is zero. A matrix with trace $T_t$ as in \eqref{eq:trace-target} has eigenvalues $e^{\lambda t/2}$ and $e^{-\lambda t/2}$, giving translation length $\lambda t$. Consequently, for every $\lambda\in[0,\sqrt2H)$, the corresponding family satisfies $t^{-1}\ell_{\Hol_X(\varphi_t(a_t))}(\gamma)=\lambda$.

By Masur's theorem \cite[Theorem~2]{Mas86}, we may choose a nonsingular, nonhorizontal periodic $q$-geodesic $\delta$ with transverse height $h>0$. The asymptotic formula \eqref{eq:periodic} then gives
\begin{equation}\label{eq:positive-anchor}
\ell_{\Hol_X(\varphi_t(a_t))}(\delta)/t\longrightarrow\sqrt2h.
\end{equation}
Since $\varphi_t(a_t)/t^2\to q$, the normalized differentials remain in a bounded subset of $Q(X)$. Lemma~\ref{lem:balanced} therefore bounds $\ell_{\Hol_X(\varphi_t(a_t))}(\beta)/t$ for every $\beta\in\Gamma$. As $\Gamma$ is countable, a diagonal argument gives a subsequence along which all normalized marked lengths converge. The resulting length function is nonzero by \eqref{eq:positive-anchor}, and \eqref{eq:trace-length} identifies its projective class with a point in the Morgan--Shalen boundary. Finally, distinct values of $\lambda$ give distinct projective classes, since the corresponding limits assign the same positive length to $\delta$ but different lengths to $\gamma$.

It remains to remove the assumption that every saddle edge of $\gamma$ is nonhorizontal. Recall that $\gamma$ was chosen using Lemma~\ref{lem:prescribed-turn} at a zero of maximal order, with the selected rays lying strictly inside two nonadjacent sectors. For small real $\theta$, set $q_\theta=e^{i\theta}q$. Since $|q_\theta|=|q|$, the curve $\gamma$ remains the same closed geodesic. In a sufficiently small interval around $\theta=0$, only finitely many values of $\theta$ make an edge of $\gamma$ horizontal. Moreover, the selected turn continues to lie in nonadjacent sectors for all sufficiently small $\theta$. Thus the preceding construction applies for every sufficiently small $\theta$ outside this finite set.

If $v_e$ is the developed flat vector of an edge, its transverse height with respect to $q_\theta$ is $|\Im(e^{i\theta/2}v_e)|$. Hence
\[
H_\theta=\sum_e|\Im(e^{i\theta/2}v_e)|\longrightarrow H>0.
\]
The positivity of $H$ already follows from the two nonhorizontal rays at the selected turn. The periodic geodesic $\delta$ remains nonsingular, and its transverse height satisfies $h_\theta\to h>0$. Given $\lambda<\sqrt2H$, we choose phases $\theta_k\to0$ for which the preceding construction applies and $\lambda<\sqrt2H_{\theta_k}$. For each $k$, this construction gives a marked length function $\ell_k$ on $\Gamma$ satisfying
\[
\ell_k(\gamma)=\lambda,\quad \ell_k(\delta)=\sqrt2h_{\theta_k}.
\]
The bounds for each element of $\Gamma$ can be chosen uniformly in $k$. Indeed, sufficiently far along each sequence realizing $\ell_k$, the differentials $\varphi/t^2$ lie in a common bounded neighborhood of $q$, and the estimate from Lemma~\ref{lem:balanced} depends only on this neighborhood and the chosen element of $\Gamma$. Since $\Gamma$ is countable, we may therefore pass to a subsequence on which $\ell_k$ converges pointwise.

Enumerate the elements of $\Gamma$ as $\{\beta_1,\beta_2,\ldots\}$. For each $k$, choose an oper from a sequence realizing $\ell_k$, with parameter $t_k\ge k$, such that
\begin{equation}\label{eq:phase-diagonal}
\|\varphi_k/t_k^2-q_{\theta_k}\|_1<1/k,\quad \left|\ell_{\Hol_X(\varphi_k)}(\beta_i)/t_k-\ell_k(\beta_i)\right|<1/k\quad(1\le i\le k).
\end{equation}
It follows that $\varphi_k/t_k^2\to q$ and the normalized marked lengths converge pointwise on $\Gamma$. The limiting length function takes the value $\lambda$ on $\gamma$ and $\sqrt2h$ on $\delta$. Finally, we choose $\lambda_k\uparrow\sqrt2H$ and repeat the diagonal selection, using the same bounds for each element of $\Gamma$ and the positivity of the $\delta$ coordinate. This also realizes $\lambda=\sqrt2H$ and completes the proof of Theorem~\ref{thm:intro-nonunique}.

\begin{proof}[Proof of Corollary~\ref{cor:intro-dichotomy}]
If all zeros of $q$ are simple, the singleton assertion follows from Propositions~\ref{prop:fiber-topology} and~\ref{prop:dumas-limits}. If $q$ has a multiple zero, Theorem~\ref{thm:intro-nonunique} gives continuum many distinct projective classes. On the other hand, every marked length function is determined by its values on the countable group $\Gamma$, so there are at most continuum many such classes. This proves the corollary.
\end{proof}


\section{Analytic realization of the boundary}\label{sec:realization}

We use a reference oper connection $\nabla_0$ on the jet bundle $E$ introduced in \eqref{eq:jet}.
In classical oper theory, the relevant criterion is whether the holomorphic bundle associated with a flat connection is isomorphic to $E$ \cite{Wen16,BHS23}. We formulate an asymptotic version of this criterion using Robinson's coefficients and a sheaf of holomorphic functions with controlled growth. The construction of the coefficient field comes from nonstandard analysis \cite{Rob66,TW04}, while the framework of function algebras has precedents in \cite{Tod99,KK18}. The passage from generalized holomorphic functions to holomorphic representatives has also been studied in \cite{OPV07}. In our construction, we perform the global correction on $E$ and estimate its effect on all monodromy matrices at the same exponential scale.


\subsection{Asymptotic coefficients and the bundle condition}\label{subsec:asymptotic-coefficients}

\begin{definition}\label{def:robinson}
Let $\tau=(\tau_j)_{j\ge1}$ be a sequence of positive real numbers with $\tau_j\to\infty$, and let $\U$ be a nonprincipal ultrafilter on $\N$. We define
\begin{align*}
\mathcal M_\tau&=\{(a_j)\in\C^{\N}: |a_j|\le e^{C\tau_j}\text{ on a set in }\U\text{ for some }C>0\},\\
\mathcal N_\tau&=\{(a_j)\in\C^{\N}: \text{for every }A>0,\ |a_j|\le e^{-A\tau_j}\text{ on a set in }\U\}.
\end{align*}
The \emph{Robinson asymptotic field} associated with $(\tau,\U)$ is the quotient $\K_{\tau,\U}=\mathcal M_\tau/\mathcal N_\tau$, equipped with the valuation
\begin{equation}\label{eq:valuation}
v([a_j])=-\lim_\U\frac{\log|a_j|}{\tau_j}\quad([a_j]\ne0),\quad v(0)=+\infty.
\end{equation}
Sequences in $\mathcal M_\tau$ are called \emph{moderate}, while those in $\mathcal N_\tau$ are called \emph{negligible}.
\end{definition}

This is the usual Robinson construction, with $[e^{-\tau_j}]$ as the infinitesimal; see \cite{Rob66,TW04}. The ideal $\mathcal N_\tau$ is maximal in $\mathcal M_\tau$, and the resulting field $\K_{\tau,\U}$ is algebraically closed. Its value group is $\R$, and its valuation is trivial on $\C$. Since a nonzero class is not negligible, every representative has an exponential lower bound on a set in $\U$, and its reciprocal is therefore moderate on that set. When $\tau$ and $\U$ are understood, we write $\K=\K_{\tau,\U}$.

\begin{definition}\label{def:holomorphic-sheaf}
For an open set $V\subset X$, let $\mathcal M_{\mathrm{hol}}(V)$ be the algebra of sequences $f_j\in\mathcal O_X(V)$ such that every compact $C^k$ seminorm is bounded by $e^{C\tau_j}$ on a set in $\U$, where $C>0$ may depend on the compact set and $k$. Let $\mathcal N_{\mathrm{hol}}(V)$ be the ideal of sequences for which, for every $A>0$, each such seminorm is bounded by $e^{-A\tau_j}$ on a set in $\U$. We define $\A_{\tau,\U}$ as the sheafification of the presheaf
\[
V\longmapsto\mathcal M_{\mathrm{hol}}(V)/\mathcal N_{\mathrm{hol}}(V).
\]
We call $\A_{\tau,\U}$ the \emph{sheaf of asymptotic holomorphic functions} associated with $(\tau,\U)$.
\end{definition}

The sets of indices may depend on the seminorm and the required rate. These conditions are preserved under coordinate changes on compact sets. Both $\mathcal O_X$ and the constant field $\K$ embed naturally in $\A_{\tau,\U}$, and differentiation is well defined because it preserves moderation and negligibility. Sheafification allows us to work with local representatives after refining the cover. Since $X$ is compact, the arguments below use finitely many relatively compact charts, and each estimate requires intersecting only finitely many sets in $\U$.

Let $\rho:\Gamma\to\SL_2(\K)$ be a representation, and let $\mathcal F_\rho$ denote the sheaf obtained by extending its local system of $\K$-vector spaces to $\A_{\tau,\U}$. This sheaf is locally free of rank two and has a determinant trivialization. On a good cover, its transition matrices $c_{ik}\in\SL_2(\K)$ are locally constant. We also write $E_\A=E\otimes_{\mathcal O_X}\A_{\tau,\U}$.

The valuation determines the standard Bruhat--Tits tree for $\SL_2(\K)$, and $\rho$ induces an action of $\Gamma$ on this tree. When the corresponding length function is nonzero, we write $T_\rho$ for its minimal invariant subtree. We normalize the metric so that $\diag(a,a^{-1})$ has translation length $2|v(a)|$. With this normalization, the marked length function is given by
\begin{equation}\label{eq:valued-length}
\ell_\rho(\gamma)=2\max\{0,-v(\tr\rho(\gamma))\}.
\end{equation}
This is the rank-two construction underlying the tree compactification \cite{MS84,CM87}. Passing to the metric completion of a tree, when needed for a Dirichlet problem, does not change its marked translation lengths.

The following is the main result of this section, which gives a necessary and sufficient analytic criterion for a projective marked length function to belong to a given directional holonomy fiber.

\begin{theorem}\label{thm:criterion}
Fix $q\in Q(X)$ with $\|q\|_1=1$. The fiber $\mathcal L_X(q)$ consists exactly of the classes $[\ell_\rho]$ for which there exist $\tau$, $\U$, and $\rho:\Gamma\to\SL_2(\K_{\tau,\U})$ satisfying the following conditions:
\begin{enumerate}
\item The tree $T_\rho$ is nontrivial and irreducible, and there is a continuous surjective $\Gamma$-equivariant straight map $T_q\to T_\rho$.
\item There is a determinant-one isomorphism of holomorphic $\A_{\tau,\U}$-bundles $\mathcal F_\rho\simeq E_\A$.
\end{enumerate}
For each such triple, there is a sequence of opers on $X$ whose monodromy over the original field $\K_{\tau,\U}$ is conjugate to $\rho$. An ordinary subsequence, indexed by $j_k$, satisfies
\begin{equation}\label{eq:criterion-convergence}
\frac{\varphi_{j_k}}{\tau_{j_k}^2}\longrightarrow\frac q2,\quad \frac{\ell_{\Hol_X(\varphi_{j_k})}(\gamma)}{\tau_{j_k}}\longrightarrow\ell_\rho(\gamma)\quad(\gamma\in\Gamma).
\end{equation}
\end{theorem}

The second condition requires the classes of $\mathcal F_\rho$ and $E_\A$ to agree in the nonabelian \v{C}ech cohomology set $H^1(X,\SL_2(\A_{\tau,\U}))$. It thus imposes a global analytic condition relating the locally constant cocycle of $\rho$ to the holomorphic cocycle of $E$. In the first condition, the metric is normalized by the source tree $T_q$, which accounts for the factor $1/2$ in \eqref{eq:criterion-convergence}.


\subsection{Construction of oper connections}\label{subsec:exact-connections}

We first use the bundle condition to construct the oper connections. The direction of the limiting quadratic differential will be recovered later.

\begin{proposition}\label{prop:controlled-descent}
Suppose that $\rho:\Gamma\to\SL_2(\K_{\tau,\U})$ is a representation and that there is a determinant-one isomorphism $\mathcal F_\rho\simeq E_\A$. Then there are tracefree holomorphic connections $\nabla_j$ on $E$ defining opers on $X$ such that their monodromy matrices are moderate and, over $\K_{\tau,\U}$, define a representation conjugate to $\rho$. The error in local transport relative to the prescribed asymptotic frames is negligible along any given compact path.
\end{proposition}

\begin{proof}
Choose a finite coordinate cover that trivializes both $E$ and $\mathcal F_\rho$. Let $g_{ik}$ be the determinant-one transition matrices of $E$, and let $c_{ik}$ be the locally constant transition matrices of the local system. After refining the cover and shrinking its members relatively compactly, we represent the given isomorphism and its inverse on the larger chart domains by holomorphic matrices $G_{i,j}$ and $H_{i,j}$ with moderate entries and derivatives. We also choose moderate constant representatives $C_{ik,j}$ of $c_{ik}$. On the compact overlaps used below, the bundle isomorphism gives
\begin{equation}\label{eq:cocycle-error}
G_{i,j}C_{ik,j}-g_{ik}G_{k,j}=R_{ik,j},
\end{equation}
with $R_{ik,j}$ negligible in every fixed $C^k$ norm. To justify these estimates after sheafification, we observe that each equality has representatives on a neighborhood of every point of an overlap. Since the overlaps are compact, finitely many such neighborhoods suffice to cover each one. Thus, for any prescribed derivative order and decay rate, we need only finitely many estimates and finitely many sets in $\U$.

The identity $G_{i,j}H_{i,j}=I+(\text{negligible})$ implies that $G_{i,j}$ is invertible on each shrunken chart for indices in a set in $\U$. The Neumann series for the error matrix shows that its inverse is moderate, and repeated differentiation of $G^{-1}$ gives moderate bounds for derivatives of every order. We also have $\det G_{i,j}=1+(\text{negligible})$. Its values lie in a disk about $1$ where the principal holomorphic inverse square root is defined. Replacing $G_{i,j}$ by $(\det G_{i,j})^{-1/2}G_{i,j}$ makes its determinant exactly one while changing it only negligibly, so \eqref{eq:cocycle-error} continues to hold. Similarly, $\det C_{ik,j}=1+(\text{negligible})$, and the adjugate formula gives moderate inverses for the constant representatives on a common set in $\U$. We choose such a set on which all these operations are defined. Outside this set, we will take the final connection to be $\nabla_0$.

We now define the local tracefree holomorphic connection forms
\begin{equation}\label{eq:pure-gauge}
A_{i,j}=-dG_{i,j}G_{i,j}^{-1}.
\end{equation}
For points $x$ and $y$ in a chart, the corresponding transport is exactly
\begin{equation}\label{eq:pure-transport}
U_{i,j}(y,x)=G_{i,j}(y)G_{i,j}(x)^{-1}.
\end{equation}
Thus both the transport and its inverse are moderate uniformly along all subpaths in any chosen compact set.

Differentiating \eqref{eq:cocycle-error}, we obtain the usual transformation law for connection forms, up to a negligible error:
\begin{equation}\label{eq:connection-overlap}
A_{i,j}=g_{ik}A_{k,j}g_{ik}^{-1}-dg_{ik}g_{ik}^{-1}+(\text{negligible}).
\end{equation}
To verify the error and its sign, set $Q=G_{i,j}C_{ik,j}$ and $S=g_{ik}G_{k,j}$. Since $C_{ik,j}$ is constant, we have $-dQ Q^{-1}=A_{i,j}$. The difference $Q-S$ and all its derivatives are negligible, whereas $Q^{-1}$ and $S^{-1}$ are moderate. Using the identity
\[
Q^{-1}-S^{-1}=-Q^{-1}(Q-S)S^{-1},
\]
we see that $-dQ Q^{-1}+dS S^{-1}$ is negligible in every prescribed seminorm. Expanding $-dS S^{-1}$ then gives \eqref{eq:connection-overlap} with the stated sign. 

Let $A_{0,i}$ be the local forms of $\nabla_0$. Then $B_{i,j}=A_{i,j}-A_{0,i}$ is a local holomorphic section of the bundle
\[
W=\End_0(E)\otimes K_X,
\]
and the differences between these sections on overlaps are negligible after applying the transition maps of $W$. We choose a smooth partition of unity $(\chi_i)$ subordinate to the shrunken cover and form a global smooth section $B_j=\sum_i\chi_i B_{i,j}$, using the appropriate transition maps in each local expression. In chart $k$, we have
\begin{equation}\label{eq:tensor-patch}
B_j-B_{k,j}=\sum_i\chi_i(B_{i,j}-B_{k,j}),\quad \dbar_W B_j=\sum_i(\dbar\chi_i)(B_{i,j}-B_{k,j}).
\end{equation}
Both expressions are negligible in every prescribed norm, since the transition maps and the derivatives of the partition of unity are independent of $j$. Moreover, $B_j$ itself is moderate.

Let $D=\dbar_W$, and let $\Pi$ be the $L^2$ orthogonal projection onto $H^0(X,W)$. Standard elliptic theory gives
\begin{equation}\label{eq:elliptic-projection}
\|B-\Pi B\|_{H^{s+1}}\le C_s\|DB\|_{H^s}\quad(s\ge0).
\end{equation}
Indeed, let $\mathcal G$ be the Green operator of $D^*D$ on the orthogonal complement of its kernel. Then $(I-\Pi)B=\mathcal G D^*DB$, and the estimate follows from elliptic regularity \cite{Wel08}. The constants depend on $W$ and its metrics, which are independent of $j$. Combining Sobolev embedding with \eqref{eq:tensor-patch}, we conclude that $B_j-\Pi B_j$ is negligible in every $C^k$ norm.

We define
\begin{equation}\label{eq:corrected-connection}
\nabla_j=\nabla_0+\Pi B_j.
\end{equation}
This is a tracefree holomorphic connection on $E$ and hence is flat. By \eqref{eq:unchanged-oper}, every global holomorphic correction tensor preserves $L$ and leaves its second fundamental form unchanged. Therefore $(E,L,\nabla_j)$ is an oper on $X$. In a local flat frame, the line bundle $L$ determines a holomorphic map to $\PP^1$ whose derivative is precisely the second fundamental form. Since this form is an isomorphism, the map is locally injective. This gives the corresponding projective structure directly.

It remains to compare the monodromy. Along any chosen chart path, the corrected connection differs from $A_{i,j}$ by a negligible form $E_j$. Let $U_j$ be the pure-gauge transport \eqref{eq:pure-transport} along this path, and let $T_j$ be the corrected transport. Writing $T_j=U_jV_j$, we obtain
\begin{equation}\label{eq:relative-transport}
V_j'=-U_j^{-1}E_jU_jV_j,\quad V_j(0)=I.
\end{equation}
For any desired decay rate, the exponential bounds on $U_j$ and $U_j^{-1}$ can be absorbed by choosing a stronger negligible bound on $E_j$. Thus the coefficient in \eqref{eq:relative-transport} is negligible. The integral inequality
\[
\|V_j-I\|\le\exp\left(\int\|U_j^{-1}E_jU_j\|\right)-1
\]
then shows that $T_j-U_j$ is negligible. This estimate uses the explicit pure-gauge transports, so its constants depend on their exponential growth rather than on the size of the untransformed connection coefficients.

Every generator of $\Gamma$ can be represented by a finite concatenation of these chart paths. Replacing the corrected transfers by the corresponding pure-gauge transfers changes each monodromy matrix negligibly. Over $\K$, equation \eqref{eq:cocycle-error} then allows the product to telescope to the prescribed locally constant cocycle, conjugated by the basepoint matrix $[G_{0,j}]$. Hence all generator matrices agree with $\rho$ up to the same conjugation over $\K$. Their entries and inverse entries are moderate. Since every word is a finite product of generator matrices and their inverses, its matrix, and therefore its trace, agrees exactly in $\K$ as well. Thus the monodromy of the constructed opers is conjugate to $\rho$ over $\K$. This completes the proof.
\end{proof}

The construction also preserves cancellation between exponentially large terms in a trace, since equality in $\K$ means that the error decays faster than every exponential rate. Applying \eqref{eq:trace-length} to the constructed connections, we obtain
\begin{equation}\label{eq:ultralimit-length}
\lim_\U\frac{\ell_{\Hol_X(\varphi_j)}(\gamma)}{\tau_j}=2\max\{0,-v(\tr\rho(\gamma))\}=\ell_\rho(\gamma).
\end{equation}
If $\tr\rho(\gamma)=0$ in $\K$, its representatives are negligible, and the limit in \eqref{eq:ultralimit-length} is zero. To obtain ordinary convergence, we enumerate $\Gamma$. For each $k$, we intersect the sets in $\U$ on which the errors in the first $k$ normalized lengths are all less than $1/k$. Since $\U$ is nonprincipal, each intersection is infinite, allowing us to choose increasing indices $j_k$. Along this subsequence, the normalized lengths converge for every element of $\Gamma$. The conjugacy of the representations over $\K$ holds on the original indexing sequence, while the chosen subsequence gives ordinary convergence of the normalized marked lengths.


\subsection{Proof of Theorem~\ref{thm:criterion}}

We now assume the first condition of Theorem~\ref{thm:criterion}: the tree $T_\rho$ is nontrivial and irreducible, and there is a straight map $T_q\to T_\rho$, where $\|q\|_1=1$. By Proposition~\ref{prop:controlled-descent} and the argument following \eqref{eq:ultralimit-length}, we may choose a subsequence of the constructed opers whose normalized marked lengths converge to $\ell_\rho$ on all of $\Gamma$. We relabel this subsequence by $j$ and set $s_j=\sqrt{\|\varphi_j\|_1}$.

\begin{lemma}\label{lem:direction-recovery}
For the subsequence chosen above, we have $\varphi_j/\tau_j^2\to q/2$.
\end{lemma}

\begin{proof}
Since $\ell_\rho$ is nonzero, there is $\gamma\in\Gamma$ for which $\ell_{\Hol_X(\varphi_j)}(\gamma)/\tau_j$ has a positive limit. Lemma~\ref{lem:balanced} gives $\ell_{\Hol_X(\varphi_j)}(\gamma)\le C_\gamma(1+s_j)$, so $s_j/\tau_j$ is bounded away from zero and $s_j\to\infty$. Suppose that $s_j/\tau_j$ were unbounded above. Passing to a subsequence, we may assume that $s_j/\tau_j\to\infty$ and $\varphi_j/\|\varphi_j\|_1\to q'$, where $\|q'\|_1=1$. Choose a nonsingular, nonhorizontal periodic $q'$-geodesic $\delta$. By \eqref{eq:periodic}, the ratio $\ell_{\Hol_X(\varphi_j)}(\delta)/s_j$ converges to a positive number. Consequently, $\ell_{\Hol_X(\varphi_j)}(\delta)/\tau_j$ would tend to infinity, contradicting the convergence of the $\tau_j$-normalized length of $\delta$. Thus
\begin{equation}\label{eq:comparable-scales}
0<c\le s_j/\tau_j\le C<\infty.
\end{equation}

Consider any subsequence along which both $\varphi_j/\|\varphi_j\|_1$ and $s_j/\tau_j$ converge, and denote the limit of the normalized differentials by $q'$. 
Applying Proposition~\ref{prop:dumas-limits} and rescaling the metrics on both trees, we obtain a straight map from $T_{q'}$ to a tree representing the same projective marked length function as $\rho$. By Lemma~\ref{lem:straight-harmonic}, the associated harmonic map has Hopf differential in the positive real ray of $-q'/4$, whereas the given straight map from $T_q$ has Hopf differential $-q/4$. The uniqueness of the normalized Hopf differential for a nonabelian projective length function \cite[Propositions~3.1 and~3.5]{DDW00}, together with the identification of minimal trees \cite[Theorem~3.7]{CM87}, therefore gives $q'=q$. Hence the normalized differentials converge to $q$.

By \eqref{eq:comparable-scales}, every accumulation point of $\varphi_j/\tau_j^2$ has the form $rq$ for some $r>0$. Let $\delta$ be a nonsingular, nonhorizontal periodic $q$-geodesic with transverse height $h>0$. Its lift projects to a complete axis in $T_q$. The straight map is isometric on each finite segment of this axis, so its image is an axis in $T_\rho$ with $\ell_\rho(\delta)=h$. On the other hand, applying \eqref{eq:periodic} along a subsequence for which $\varphi_j/\tau_j^2\to rq$, we obtain
\[
\ell_{\Hol_X(\varphi_j)}(\delta)/\tau_j\longrightarrow\sqrt{2r}\,h.
\]
Comparing the two limits gives $r=1/2$. Thus every accumulation point of $\varphi_j/\tau_j^2$ equals $q/2$, which proves the lemma.
\end{proof}

\begin{proof}[Proof of Theorem~\ref{thm:criterion}]
We first prove sufficiency. Under the second condition, Proposition~\ref{prop:controlled-descent} constructs a sequence of opers whose monodromy over $\K_{\tau,\U}$ is conjugate to $\rho$. As shown after \eqref{eq:ultralimit-length}, we may pass to an ordinary subsequence along which the normalized marked lengths converge to $\ell_\rho$ for every element of $\Gamma$. The first condition ensures that this limit is nonzero, while Lemma~\ref{lem:direction-recovery} gives \eqref{eq:criterion-convergence}. Since $\|q\|_1=1$, the Schwarzians diverge and approach the positive ray of $q$. Hence $[\ell_\rho]\in\mathcal L_X(q)$.

For necessity, let $[\ell]\in\mathcal L_X(q)$ and choose a sequence $\varphi_j\in Q(X)$ realizing this boundary point. Set $\tau_j=\sqrt{2\|\varphi_j\|_1}$. The upper bounds on the lengths and the countability of $\Gamma$ allow us to pass to a subsequence along which all $\tau_j$-normalized marked lengths converge. The asymptotic formula for a periodic $q$-geodesic gives a positive limiting length, and \eqref{eq:scale-comparison} shows that the resulting length function represents $[\ell]$. We also have $\varphi_j/\tau_j^2\to q/2$.

Choose a nonprincipal ultrafilter $\U$ on this indexing set and realize the scalar equations on the jet bundle $E$. By Lemma~\ref{lem:balanced}, the exact holomorphic solution matrices and their inverses are bounded by $e^{C\tau_j}$ on a finite collection of coordinate charts and connecting paths. After shrinking the disks, Cauchy's estimates give the same exponential bounds for holomorphic derivatives of every prescribed order. Choosing initial conditions of determinant one gives solution frames with determinant one. The exact monodromy identities pass to $\K_{\tau,\U}$ and define a representation $\rho:\Gamma\to\SL_2(\K_{\tau,\U})$. Moreover, the local solution frames give an isomorphism $\mathcal F_\rho\simeq E_\A$ of determinant one. This establishes the second condition.

Since $\varphi_j/\tau_j^2\to q/2$, Proposition~\ref{prop:dumas-limits} gives a surjective $\Gamma$-equivariant straight map from $T_q$ to a tree representing the limiting marked lengths. By Lemma~\ref{lem:straight-harmonic}, the composite map from $\widetilde X$ is harmonic, and its image is minimal by the remark following \cite[Proposition~2.4]{DDW00}. Corollary~\ref{cor:irreducible-limits} shows that the action preserves neither an end nor a line, even after restriction to a finite-index subgroup. By \eqref{eq:valued-length} and \eqref{eq:trace-length}, the valued tree $T_\rho$ has the same marked length function. The rigidity theorem for irreducible actions \cite[Theorem~3.7]{CM87} therefore gives an equivariant identification of the two minimal trees, establishing the first condition. Finally, all estimates in the necessity argument concern exact scalar solutions on coordinate disks chosen independently of $j$, including disks containing zeros. Thus the argument applies to every sequence realizing a point of $\mathcal L_X(q)$, completing the proof.
\end{proof}

The following proposition relates the bundle condition to ordinary holomorphic vector bundles.

\begin{proposition}\label{prop:fixed-bundles}
If two complex holomorphic vector bundles $F$ and $G$ on $X$ become isomorphic after extension to $\A_{\tau,\U}$, then $F\simeq G$ as holomorphic vector bundles on $X$.
\end{proposition}

\begin{proof}
On a finite cover with relatively compact members, we represent the isomorphism and its inverse by holomorphic matrix maps with moderate entries and derivatives. Their differences on overlaps are negligible. Using a partition of unity, as in \eqref{eq:tensor-patch}, we obtain a global smooth map $H_j:F\to G$ whose Dolbeault derivative is negligible. Since $H_j$ differs negligibly from each invertible local map, it is invertible for indices in a set in $\U$, and its inverse is moderate. We then apply the Dolbeault projection on $\Hom(F,G)$. By \eqref{eq:elliptic-projection}, the resulting global holomorphic map $H_j^{\mathrm{hol}}$ differs negligibly from $H_j$. Consequently,
\[
\|H_j^{-1}(H_j^{\mathrm{hol}}-H_j)\|<1
\]
uniformly on $X$ for indices in a set in $\U$. For any such index, $H_j^{\mathrm{hol}}$ is a holomorphic bundle isomorphism, proving the proposition.
\end{proof}

Consequently, for a constant complex representation, the second condition of Theorem~\ref{thm:criterion} reduces to the classical criterion for an oper in terms of an isomorphism of holomorphic bundles. For a representation over asymptotic coefficients, the same condition instead controls holomorphic representatives of its cocycle at every exponential rate. The construction above uses these representatives to obtain exact oper connections and determine the limiting length function on every element of $\Gamma$.

\vspace{3mm}

\noindent 
\textbf{Acknowledgment.} S.T. is supported by NSFC (Grant No.~12601116) and BJNSF (Grant No.~1264052).

\vspace{2mm}

\noindent\textbf{AI disclosure.}
The authors used OpenAI ChatGPT (GPT-6) to explore proof strategies, critically examine mathematical arguments, and assist with revisions of the first manuscript they wrote. All mathematical arguments and proofs were independently verified by the authors, who take full responsibility for the mathematical results and the final text.

\end{document}